\documentclass{compositio}
\usepackage{amssymb,latexsym,amsmath,amscd}

\theoremstyle{theorem}
\newtheorem{theorem}{Theorem}[section]
\newtheorem{corollary}{Corollary}[section]
\newtheorem{lemma}{Lemma}[section]
\newtheorem{proposition}{Proposition}[section]

\theoremstyle{definition}
\newtheorem{definition}{Definition}[section]
\newtheorem{example}{Example}[section]

\theoremstyle{remark}
\newtheorem{remark}{Remark}[section]
\newtheorem{notation}{Notation}[section]

\numberwithin{equation}{subsection}

\newcommand{\cA}{{\mathcal A}}
\newcommand{\cC}{{\mathcal C}}
\newcommand{\cE}{{\mathcal E}}
\newcommand{\cF}{{\mathcal F}}
\newcommand{\cO}{{\mathcal O}}
\newcommand{\cP}{{\mathcal P}}
\newcommand{\cQ}{{\mathcal Q}}

\newcommand{\cT}{{\mathcal T}}
\newcommand{\cU}{{\mathcal U}}
\newcommand{\cV}{{\mathcal V}}

\newcommand{\im}{\operatorname{Im}}

\newcommand{\Diff}{\operatorname{Diff}}

\newcommand{\isomo}{\overset{\sim}{=}}

\newcommand{\id}{\operatorname{id}}
\newcommand{\Tr}{\operatorname{Tr}}
\newcommand{\shHom}{\underline{\operatorname{Hom}}}
\newcommand{\shEnd}{\underline{\operatorname{End}}}
\newcommand{\shExt}{\underline{\operatorname{Ext}}}
\newcommand{\shAut}{\underline{\operatorname{Aut}}}

\newcommand{\ch}{\operatorname{ch}}

\newcommand{\pr}{\operatorname{pr}}

\newcommand{\ECA}{\mathcal{ECA}}
\newcommand{\VA}{\mathcal{VA}}
\newcommand{\EVA}{\mathcal{EVA}}
\newcommand{\CEXT}{\mathcal{CE\!\scriptstyle{XT}}}
\newcommand{\VEXT}{\mathcal{VE\!\scriptstyle{XT}}}
\newcommand{\EXT}{\mathcal{E\!\scriptstyle{XT}}}
\newcommand{\DR}{{\Omega^\bullet_X}}

\newcommand{\vac}{{\mathbf{1}}}
\newcommand\dual[1]{{#1}^{\vee}}
\newcommand{\ip}{{\langle\ ,\ \rangle}}

\newcommand{\pont}{{\scriptstyle{\Pi}}}

\begin{document}

\title{The first Pontryagin class}

\author[P. Bressler]{Paul Bressler}
\email{bressler@math.ias.edu}
\address{School of Mathematics\\
         Institute for Advanced Study\\
         Princeton, NJ 08540}

\begin{abstract}
We give a natural obstruction theoretic interpretation to the first Pontryagin
class in terms of Courant algebroids. As an application we calculate the class of
the stack of algebras of chiral differential operators. In particular, we establish
the existence and uniqueness of the chiral de Rham complex.
\end{abstract}


\maketitle

\section{Introduction}
\subsection{The first Pontryagin class}
The first Pontryagin class, for the purposes of the present paper, is a characteristic class associated to a pair consisting of
a principal $G$-bundle, $G$ a Lie group, over a manifold $X$ and an invariant symmetric bilinear form $\ip$ on the Lie algebra
$\mathfrak{g}$ of $G$. For a $G$-bundle $P$ on $X$ the Pontryagin class, denoted by $\pont(P,\ip)$ takes values in
$H^2(X;\Omega^2\to\Omega^{3,cl})$.

Incarnations of the first Pontryagin class corresponding to particular choices of $(G,\ip)$ are quite familiar. For example, the class $2\ch_2$ is the Pontryagin class
corresponding to $GL_n(\mathbb{C})$ and the canonical pairing on
$\mathfrak{gl}_n$ given by the trace of the product of matrices.

More generally, the first Pontryagin class with values as above may be associated to
a transitive Lie algebroid (see \ref{ssection:lie-algebroids}), say, $\cA$, together with an invariant symmetric pairing $\ip$ on the kernel of the anchor map and will be denoted $\pont(\cA,\ip)$.\footnote{The construction of characteristic classes in this general context and the calculation of the \v{C}ech-de Rham representative of the first Pontryagin class are recalled in Appendix \ref{section:char-classes}.} The Pontryagin class of a principal bundle is defined as the Pontryagin class of the Atiyah algebra of the bundle.

\subsection{The first Pontryagin class as an obstruction}
It turns out that the Pontryagin class $\pont(\cA,\ip)$ appears naturally in the context of a certain classification problem which is canonically associated with the pair $(\cA,\ip)$.

Just as the degree one cohomology classifies torsors (principal bundles), degree two
cohomology classifies certain stacks (\cite{Br}).
To each transitive Lie algebroid $\cA$ on $X$ and an invariant symmetric bilinear form
$\ip$ on the kernel of the anchor map $\pi:\cA\to\cT_X$
we will associate the stack $\CEXT_{\cO_X}(\cA)_\ip$ with the corresponding class in
$H^2(X;\Omega^2\to\Omega^{3,cl})$ equal to $-\dfrac12\pont(\cA,\ip)$. The latter equality is the content of Theorem \ref{thm:main}.

The stack $\CEXT_{\cO_X}(\cA)_\ip$ associates to an open subset $U$ of $X$ the category (groupoid) $\CEXT_{\cO_X}(\cA)_\ip(U)$ of certain \textit{Courant extensions} of $\cA$. Courant extensions of Lie algebroids are, in particular, \textit{Courant algebroids}. Consequently, a significant part of the paper is devoted to the definition, basic properties and classification of Courant algebroids, the main result (Theorem \ref{thm:main}) being the identification of the class of $\CEXT_{\cO_X}(\cA)_\ip$ with $-\dfrac12\pont(\cA,\ip)$.

Thus, $(\cA,\ip)$ admits a (globally defined) Courant extension if and only if the the Pontryagin class of $(\cA,\ip)$ vanishes.
This fact (sans the terminology of the present paper) is indicated in \cite{S}.
\subsection{Courant extensions of Atiyah algebras}
The following is intended to convey in an informal fashion the differential geometric meaning of
a Courant extension of the Atiyah algebra of a principal bundle.

Suppose that $G$ is a reductive group with Lie algebra $\mathfrak{g}$. Let $\underline{G}$ denote the sheaf
of groups represented by $G$. An invariant symmetric pairing $\ip$ on $\mathfrak{g}$ which satisfies
certain integrality conditions gives rise to a central extension, say $\widehat{\underline{G}}_\ip$ of
$\underline{G}$ by (the sheaf) $\underline{K}_2$ so that there is a short exact sequence of sheaves of groups
\[
\begin{CD}
1 @>>> \underline{K}_2 @>>> \widehat{\underline{G}}_\ip @>>> \underline{G} @>>> 1
\end{CD}
\]
which leads to the short exact sequence of tangent spaces at the identity
\[
\begin{CD}
0 @>>> \Omega^1 @>>> \widehat{\underline{\mathfrak{g}}}_\ip @>>> \underline{\mathfrak{g}} @>>> 0
\end{CD}
\]
where $T_e \underline{K}_2 = \Omega^1$, $\widehat{\underline{\mathfrak{g}}}_\ip$
denotes $T_e\widehat{\underline{G}}_\ip$ and $\underline{\mathfrak{g}} = T_e\underline{G}$ is the sheaf of Lie
algebras represented by $\mathfrak{g}$. As was pointed out by S.~Bloch in \cite{B} (at least in the case of
the Steinberg group), the usual construction does not yield a Lie bracket on $\widehat{\underline{\mathfrak{g}}}_\ip$.
Thus, $\widehat{\underline{\mathfrak{g}}}_\ip$ is not a sheaf of Lie algebras in a way compatible with the projection
to $\underline{\mathfrak{g}}$. It \textit{is} however a Courant algebroid (with the trivial anchor map) and, in
particular, a sheaf of Leibniz algebras (the Leibniz bracket, however, is \textit{not} $\cO$-linear). Examples
of this kind are studied in \ref{sssec:LtoL}.

Now suppose that $P$ is a $\underline{G}$-torsor  (i.e. a principal $G$-bundle) on $X$. Lifts of $P$
a  $\widehat{\underline{G}}_\ip$ (i.e. pairs $(\widehat{P},\phi)$ comprised of a $\widehat{\underline{G}}_\ip$-torsor
$\widehat{P}$ and a $\widehat{\underline{G}}_\ip$-equivariant map $\phi : \widehat{P}\to P$) exist locally
on $X$ and form a $\underline{K}_2$-gerbe whose class in $H^2(X;\underline{K}_2)$ is the class of the central
extension.

Given a lift $\widehat{P}$ as above one might ask what a connection on $\widehat{P}$ might be, or, better,
what sort of structure would serve as the ``Atiyah algebra" in this context. Denoting this hypothetical, for the
moment, object by $\widehat\cA_P$ we note that it may be expected, at the very least, to fit into the commutative
diagram
\[
\begin{CD}
\Omega^1 @>>> \widehat{\underline{\mathfrak{g}}}_\ip @>>> \underline{\mathfrak{g}}\\
@V{\id}VV @VVV  @VVV \\
\Omega^1 @>>> \widehat\cA_P @>>>\cA_P \\
& & @VVV @VVV \\
& & \cT_X @>{\id}>> \cT_X
\end{CD}
\]
(where $\cA_P$ denotes the Atiyah algebra of $P$) and carry a (Leibniz) bracket compatible with all of
the maps. The above picture is encapsulated in the notion of a \textit{Courant extension} of the Lie algebroid
$\cA_P$.

Courant extensions of $\cA_P$ which fit into the above diagram exist
locally on $X$ (due to local triviality of $P$) and form the gerbe
which we denoted by  $\CEXT_{\cO_X}(\cA_P)_\ip$ above.

\subsection{Vertex algebroids}
In the rest of the paper is devoted to an application of the thus far developed theory of Courant algebroids and
classification thereof to the question of classification of \textit{exact vertex algebroids} or, equivalently,
sheaves of chiral differential operators. Vertex algebroids were defined and classified in terms of local data
in \cite{GMS} (and other papers by the same authors) and, later, in the language of chiral algebras, in \cite{BD}.

The approach to the classification of vertex algebroids carried out in the present paper is suggested by some
results contained in \cite{BD}\footnote{The author is grateful to A.~Beilinson for sending him an early preprint
of ``Chiral algebras".}.  It turns out that the latter problem reduces (in the sense of equivalence of stacks,
Proposition \ref{prop:CEXToppEVA}) to
the classification problem for Courant extensions of $\cA_{\Omega^1_X})_\ip$ (the Atiyah algebra of the cotangent
bundle with  the symmetric pairing $\ip$ on the Lie algebra $\shEnd_{\cO_X}(\Omega^1_X)$ given by the trace
of the product of endomorphisms). The result (Theorem \ref{main-thm}) is that the the stack $\EVA_{\cO_X}$ of
vertex algebroids on $X$ gives rise to a class (obstruction to existence of a globally defined vertex algebroid)
in $H^2(X;\Omega^2\to\Omega^{3,cl})$, and that class is equal to $\ch_2(\Omega^1_X)$.

This reduction is achieved with the aid of (the degree zero component of) the differential graded vertex algebroid
over the de Rham complex of $X$ (which gives rise to the chiral de Rham complex of \cite{MSV}). The existence
and uniqueness of this object was demonstrated in \cite{MSV} (and \cite{BD} in the language of chiral algebras).
We give a ``coordinate-free'' proof of this result. To this end we apply the results of the preceding sections to
the differential graded manifold $X^\sharp$ whose underlying space is $X$ and the structure sheaf is the de Rham
complex of $X$. We show (Proposition \ref{prop:ECA-DR-triv}) that every (differential graded) exact Courant algebroid
on $X^\sharp$ is canonically trivialized. This implies (Corollary \ref{corollary:existence-uniqueness}) that an exact
vertex algebroid on $X^\sharp$ exists and is unique up to a unique isomorphism. This is a differential graded object
whose degree zero constituent is a vertex extension of the Atiyah algebra of the cotangent bundle.

\subsection{Acknowledgements} The author would like to thank Pavol \v{S}evera
for many inspiring discussions. Much of the present work was carried
out during the author's visits to Universit\'{e} d'Angers and
I.H.\'E.S.

\section{Courant algebroids}

\subsection{Leibniz algebras}
Suppose that $k$ is a commutative ring.

\begin{definition}
A {\em Leibniz $k$-algebra} is a $k$-module $\mathfrak{g}$ equipped with
a bilinear operation
\[
[\ ,\ ]:\mathfrak{g}\otimes_k \mathfrak{g}\to\mathfrak{g}
\]
(the Leibniz bracket) which satisfies the Jacobi type identity
\[
[a,[b,c]] = [[a,b],c] + [b,[a,c]]\ .
\]

A morphism of Leibniz $k$-algebras is a $k$-linear map which commutes with
the respective Leibniz brackets.
\end{definition}

\begin{example}
Suppose that $\mathfrak{g}$ is a Lie algebra, $\widehat{\mathfrak{g}}$ is a
$\mathfrak{g}$-module, and $\pi:\widehat{\mathfrak{g}} \to \mathfrak{g}$ is a morphisms of
$\mathfrak{g}$-modules. The bilinear operation on $\widehat{\mathfrak{g}}$ defined by
\[
[a,b] = \pi(a)(b) \ .
\]
for $a,b\in\widehat{\mathfrak{g}}$ satisfies the Jacobi identity and thus defines
a structure of a Leibniz algebra on $\widehat{\mathfrak{g}}$.
\end{example}

\subsection{Courant algebroids}
Courant algebroids, as defined below appear as quasi-classical limits of the vertex algebroids
(see Definition \ref{definition:va} and \ref{ssection:quantCa} for discussion of quantization).
The format of the definition given below follows that of the corresponding non-commutative notion
(``vertex algebroid") which, in turn, is distilled from the structure of a vertex operator algebra
in \ref{ssection:VA-to-va}.

\begin{definition}
A {\em Courant $\cO_X$-algebroid} is an $\cO_X$-module $\cQ$
equipped with
\begin{enumerate}
\item a structure of a Leibniz $\mathbb{C}$-algebra
\[
[\ ,\ ] : \cQ\otimes_\mathbb{C}\cQ \to \cQ \ ,
\]

\item
an $\cO_X$-linear map of Leibniz algebras (the anchor map)
\[
\pi : \cQ \to \cT_X \ ,
\]

\item
a symmetric $\cO_X$-bilinear pairing
\[
\ip : \cQ\otimes_{\cO_X}\cQ \to \cO_X \ ,
\]

\item
a derivation
\[
\partial : \cO_X \to \cQ
\]
\end{enumerate}
which satisfy
\begin{eqnarray}
\pi\circ\partial & = & 0 \label{complex}\\
\left[q_1,fq_2\right] & = & f[q_1,q_2] + \pi(q_1)(f)q_2 \label{leibniz}\\
\langle [q,q_1],q_2\rangle + \langle q_1,[q,q_2]\rangle & = & \pi(q)(\langle q_1, q_2\rangle)
\label{ip-invar}\\
\left[q,\partial(f)\right] & = & \partial(\pi(q)(f)) \label{bracket-o-courant}\\
\langle q,\partial(f)\rangle & = & \pi(q)(f) \label{ip-o}\\
\left[q_1,q_2\right] + [q_2,q_1] & = & \partial(\langle q_1, q_2\rangle) \label{ip-symm}
\end{eqnarray}
for $f\in\cO_X$ and $q,q_1,q_2\in\cQ$.

A morphism of Courant $\cO_X$-algebroids is an $\cO_X$-linear map
of Leibnitz algebras which commutes with the respective anchor
maps and derivations and preserves the respective pairings.
\end{definition}

\begin{remark}
The definition of \textit{Courant algebroid} given below reduces to Definition 2.1 of \cite{LWX}
under the additional hypotheses of loc. cit. (that $\cQ$ is locally free of finite rank and the
symmetric pairing is non-degenerate).

Courant algebroids are to \textit{vertex Poisson algebras} (coisson algebras in the terminology of
\cite{BD}) what vertex algebroids are to vertex algebras in the sense of the analysis carried out at
the outset of Section \ref{section:VA}.
\end{remark}

\subsection{Twisting by $3$-forms}\label{twist}
Suppose that $\cQ$ is a Courant algebroid with Leibniz bracket denoted $[\ ,\ ]$,
and $H$ is a $3$-form on $X$. Let $[\ ,\ ]_H$ denote the bilinear operation on
$\cQ$ defined by the formula
\begin{equation}\label{formula:twisted-by-H-bracket}
[q_1,q_2]_H = [q_1,q_2] + i(\iota_{\pi(q_2)}\iota_{\pi(q_1)}H) \ ,
\end{equation}
where $i:\Omega^1_X\to\cQ$ is the $\cO_X$-linear map in the canonical factorization of
the derivation $\partial$.

Recall that the Jacobiator $J(\{\ ,\ \})$ of a binary operation $\{\ ,\ \}$
is defined by the formula
\begin{equation}\label{formula:jacobiator}
J(\{\ ,\ \})(a,b,c) =\{a,\{b,c\}\} - \{\{a,b\},c\} - \{b,\{a,c\}\} \ .
\end{equation}

\begin{lemma}
$J([\ ,\ ]_H = i(dH)$
\end{lemma}
\begin{proof}
Direct calculation
\end{proof}

\begin{notation}
For a Courant algebroid $\cQ$ and a closed $3$-form $H$ on $X$ we denote by $\cQ_H$ the Courant
algebroid with the underlying $\cO_X$-module $\cQ$ equipped with the same symmetric pairing and derivation,
and with the Leibniz bracket $[\ ,\ ]_H$ given by \eqref{formula:twisted-by-H-bracket}.

We refer to $\cQ_H$ as \textit{the $H$-twist of} $\cQ$.
\end{notation}

It is clear that twisting by a $3$-form is a functorial operation: a morphism of Courant algebroids is also a morphsim
of their respective twists (by the same form).

\subsection{The associated Lie algebroid}
Suppose that $\cQ$ is a Courant $\cO_X$-algebroid.
\begin{notation}
Let $\Omega_\cQ$ denote the $\cO_X$-submodule of $\cQ$ generated by the image
of the derivation $\partial$. Let $\overline\cQ = \cQ/\Omega_\cQ$.
\end{notation}

For $q\in\cQ$, $f,g\in\cO_X$
\begin{eqnarray*}
[q,f\partial(g)] & = & f[q,\partial(g)] + \pi(q)(f)\partial(g) \\
& = & f\partial(\pi(q)(g)) + \pi(q)(f)\partial(g)
\end{eqnarray*}
which shows that $[\cQ,\Omega_\cQ]\subseteq\Omega_\cQ$. Therefore,
the Leibniz bracket on $\cQ$ descends to the bilinear operation
\begin{equation}\label{LAbracket}
[\ ,\ ] : \overline\cQ\otimes_\mathbb{C}\overline\cQ \to \overline\cQ\ .
\end{equation}

Since $\pi$ is $\cO_X$-linear and $\pi\circ\partial = 0$, $\pi$
vanishes on $\Omega_\cQ$, hence, factors through the map
\begin{equation}\label{LAanchor}
\pi : \overline\cQ \to \cT_X \ .
\end{equation}

\begin{lemma}
The bracket \eqref{LAbracket} and the anchor \eqref{LAanchor}
determine a structure of a Lie $\cO_X$-algebroid on
$\overline\cQ$.
\end{lemma}
\begin{proof}
According to \eqref{ip-symm} the symmetrization of the Leibniz
bracket on $\cQ$ takes values in $\Omega_{\cQ}$. Therefore, the
induced bracket is skew-symmetric. The Leibniz rule and the Jacobi
identity for $\overline\cQ$ follow from those for $\cQ$.
\end{proof}

In what follows we refer to the Lie algebroid $\overline\cQ$ as
{\em the Lie algebroid associated to the Courant algebroid $\cQ$}.

\begin{definition}
A Courant extension of a Lie algebroid $\cA$ is a Courant algebroid
$\cQ$ together with an isomorphism
$\overline{\cQ} =\cA$ of Lie algebroids.

A morphism $\phi : \cQ_1\to\cQ_2$ of Courant extensions of
$\cA$ is a morphism of Courant algebroids which is compatible with the
identifications $\overline{\cQ_i}=\cA$.
\end{definition}

Suppose that $\cA$ is a Lie algebroid. For each open subset $U$ of $X$
there is category of Courant extensions of $\cA\vert_U$. Together with
the obvious restriction functors these form a stack.
\begin{notation}
We denote the stack of Courant extensions of $\cA$ by
$\CEXT_{\cO_X}(\cA)$.
\end{notation}

\subsection{From Leibniz to Lie}
For a Lie algebroid $\cA$ (respectively, Courant algebroid $\cQ$) we will denote
by ${\mathfrak{g}}(\cA)$ (respectively, ${\mathfrak{g}}(\cQ)$) the kernel of the anchor map of $\cA$ (respectively, $\cQ$). Note that ${\mathfrak{g}}(\cQ)$ is,
naturally, a Courant algebroid with the trivial anchor map and
$\overline{{\mathfrak{g}}(\cQ)} = {\mathfrak{g}}(\overline{\cQ})$. Since
$\langle{\mathfrak{g}}(\cQ),\Omega_\cQ\rangle = 0$, the pairing $\ip$
on $\cQ$ induces the pairings
\begin{eqnarray}
\ip : {\mathfrak{g}}(\cQ)\otimes_{\cO_X}\overline{\cQ} & \to & \cO_X
\label{pairing:gA}
\\
\ip : \mathfrak{g}({\overline{\cQ}})\otimes_{\cO_X}\mathfrak{g}({\overline{\cQ}})
& \to & \cO_X \label{pairing:gg}\ .
\end{eqnarray}

The (restriction of the left) adjoint action
$\cQ\to\shEnd_\mathbb{C}(\mathfrak{g}(\cQ))$ is a
morphism of Leibniz algebras (this is equivalent to the Jacobi identity)
which annihilates $\Omega_\cQ$, hence, factors through the morphism of Lie
algebras
\begin{equation}\label{Q-bar-action}
\overline{\cQ}\to\shEnd_\mathbb{C}(\mathfrak{g}(\cQ))
\end{equation}
or, equivalently, induces a canonical structure of a Lie $\overline{\cQ}$-module
on $\mathfrak{g}(\cQ)$.

\subsection{Transitive Courant algebroids}
\begin{definition}
A Courant $\cO_X$-algebroid is called {\em transitive} if the anchor
is surjective.
\end{definition}

\begin{remark}
A Courant $\cO_X$-algebroid is transitive if and only if the
associated Lie algebroid is.
\end{remark}

Suppose  that  $\cQ$ is a transitive Courant $\cO_X$-algebroid.

\begin{lemma}\label{lemma:exactness of the exact}
Suppose that $\cQ$ is a transitive Courant algebroid. Then, the sequence
\begin{equation}\label{ses:assoc-Lie-alg}
0 \to \Omega^1_X \to \cQ \to \overline\cQ \to 0
\end{equation}
where the first map is induced by the derivation $\partial$ and the
second map is the canonical projection is exact. Moreover, the
embedding $\Omega^1_X \to \cQ$ is isotropic with respect to the
symmetric pairing.
\end{lemma}
\begin{proof}
It suffices to note that $\pi\circ i = 0$ and check that $i$ is
a monomorphisms. Since $\langle q, i(\alpha)\rangle = \iota_{\pi(q)}\alpha$,
it follows that the map $i$ is adjoint to the anchor map $\pi$. The
surjectivity of the latter implies that $i$ is injective.
\end{proof}

In what follows, when dealing with a transitive Courant algebroid
$\mathcal{Q}$ we will, in view of Lemma \ref{lemma:exactness of the
exact}, regard $\Omega^1_X$ as a subsheaf of $\mathcal{Q}$, i.e.
identify the former with its image under the embedding into the
latter. This should not lead to confusion since morphisms of
transitive Courant algebroids will, under the above identification,
induce the identity map on $\Omega^1_X$.

\begin{remark}\label{remark:CEXT-gpd}
The exact sequence \eqref{ses:assoc-Lie-alg} is functorial. Thus,
a morphism of Courant extensions of a transitive Lie algebroid
$\cA$ induces a morphism of associated extensions of $\cA$ by
$\Omega^1_X$. A morphism of extensions is necessarily an
isomorphism on the respective middle terms, and it is clear that
the inverse isomorphism is a morphism of Courant extensions of
$\cA$. Hence, the category of Courant extensions of a transitive
Lie algebroid is, in fact, a groupoid.
\end{remark}

\begin{definition}
A {\em connection} on a transitive Courant $\cO_X$-algebroid $\cQ$
is a $\cO_X$-linear isotropic (with respect to the symmetric pairing on $\cQ$)
section of the anchor map $\cQ \to \cT_X$.

A {\em flat connection} on a transitive Courant $\cO_X$-algebroid
$\cQ$ is a $\cO_X$-linear section of the anchor map which is
morphism of Leibniz algebras.
\end{definition}

\begin{lemma}
A flat connection is a connection.
\end{lemma}
\begin{proof}
Suppose that $\nabla$ is a flat connection on $\cQ$. Then,
$\im(\nabla)$ is a \emph{Lie} subalgebra in $\cQ$, i.e. the
restriction of the bracket to $\im(\nabla)$ is skew-symmetric. It
follows from \eqref{ip-symm} that the ($\mathbb{C}$-linear)
composition $\partial\circ\ip\circ(\nabla\otimes\nabla) :
\mathcal{T}_X\otimes_\mathbb{C}\mathcal{T}_X\to\mathcal{Q}$ is equal
to zero. Therefore, the ($\mathcal{O}_X$-linear) composition
$\ip\circ(\nabla\otimes\nabla) :
\mathcal{T}_X\otimes\mathcal{T}_X\to\mathcal{O}_X$ factors through
the inclusion $\mathbb{C}\to\mathcal{O}_X$, hence is trivial.

Thus, $\ip$ vanishes on $\im(\nabla)$, i.e. the latter is isotropic,
which means that $\nabla$ is a connection.
\end{proof}

\begin{lemma}\label{lemma:loc-conn}
A transitive Courant algebroid which is a locally free $\cO_X$ module admits
a connection locally on $X$.
\end{lemma}
\begin{proof}
Suppose that $\cQ$ is a Courant algebroid as above. Let $s : \cT_X\to\cQ$ denote
a locally defined section of the anchor map (such exist in a neighborhood of every
point of $X$). Let $\phi : \cT_X\to\Omega^1_X$ be defined by
\[
\iota_\eta\phi(\xi) = - \frac12\langle s(\xi),s(\eta)\rangle\ .
\]
Then, as is easy to check, $s + \phi$ is a connection.
\end{proof}

\begin{notation}
We denote by $\cC(\cQ)$ (respectively, $\cC^\flat(\cQ)$) the sheaf of
(locally defined) connections (respectively, flat connections) on $\cQ$.
\end{notation}

\begin{definition}\label{definition:curvature courant}
For a connection $\nabla$ on $\mathcal{Q}$ the \emph{curvature of
$\nabla$}, denoted $c(\nabla)$, is defined by the formula
\[
c(\nabla)(\xi,\eta) = [\nabla(\xi),\nabla(\eta)] -
\nabla([\xi,\eta]) \ ,
\]
for $\xi,\eta\in\mathcal{T}_X$.
\end{definition}

\begin{remark}\label{remark:curvature}
The curvature of connections on Courant algebroids shares basic
properties with the corresponding notion for Lie algebroids
(\ref{connections for Lie}).
\begin{enumerate}
\item It is clear that $c(\nabla)$ takes values in
$\mathfrak{g}(\mathcal{Q})$. The usual calculations show that
$c(\nabla)$ is $\mathcal{O}_X$-bilinear. Since it is clearly
alternating, it determines a map
$c(\nabla):\bigwedge^2\mathcal{T}_X\to\mathfrak{g}(\mathcal{Q})$ or,
equivalently a section
$c(\nabla)\in\Gamma(X;\Omega^2_X\otimes\mathfrak{g}(\mathcal{Q})$.

\item A connection is flat if and only if its curvature is equal to
zero.
\item Suppose that $\nabla$ is a connection on $\mathcal{Q}$. Let
$\overline{\nabla}$ denote the composition of $\nabla$ with the
projection $\overline{(\bullet)} :
\mathcal{Q}\to\overline{\mathcal{Q}}$. Then, $\overline{\nabla}$ is
a connection on the Lie algebroid $\overline{Q}$ and
$c(\overline{\nabla}) = \overline{c(\nabla)}$.
\end{enumerate}
\end{remark}

\section{Courant extensions}
\subsection{Courant extensions of transitive Lie
algebroids}\label{subsection:CEXT}
From now on we assume that
$\cA$ is a {\em transitive} Lie $\cO_X$-algebroid locally free of
finite rank over $\cO_X$. By Remark \ref{remark:CEXT-gpd},
$\CEXT_{\cO_X}(\cA)$ is a stack in groupoids.

Suppose that $\widehat\cA$ is a Courant extension of $\cA$. Let
\begin{equation}\label{map-exp}
\exp : \Omega^2_X\to\shAut_{\shExt^1_{\cO_X}(\cA,\Omega^1_X)}(\widehat\cA)
\end{equation}
denote the map defined by $\exp(B)(a) = a +\iota_{\pi(a)}B$.

\begin{lemma}\label{lemma:exp-iso}
The map \eqref{map-exp} establishes an isomorphism
\[
\exp : \Omega^{2,cl}_X\to\shAut_{\CEXT_{\cO_X}(\cA)}(\widehat\cA) \ .
\]
\end{lemma}
\begin{proof}
Suppose that $\phi\in\shAut_{\CEXT_{\cO_X}(\cA)}(\widehat\cA)$.
Then, $\phi$ restricts to the identity on $\Omega^1$ and induces the
identity on $\cA$. Therefore, $\phi(a) = a + \phi^\prime(\pi(a))$, where
$\phi^\prime : \cT_X\to\Omega^1_X$. Since, for $a_1, a_2\in\widehat\cA$,
\begin{multline*}
\langle\phi(a_1),\phi(a_2)\rangle = \langle a_1, a_2\rangle +
\langle\phi^\prime(\pi(a_1)), a_2\rangle + \langle a_1, \phi^\prime(\pi(a_2))\rangle =
\\
\langle a_1, a_2\rangle + \iota_{\pi(a_2)}\phi^\prime(\pi(a_1)) +
\iota_{\pi(a_1)}\phi^\prime(\pi(a_2))
\end{multline*}
it follows that $\phi$ preserves the symmetric pairing if and only if $\phi^\prime$, viewed as a
section of $\Omega^1_X\otimes_{\cO_X}\Omega^1_X$, is alternating, i.e.
$\phi^\prime = B\in\Omega^2_X$ and $\phi = \exp(B)$.

The formula\footnote{computed using the identity $d\iota_\eta\iota_\xi =
(L_\eta - \iota_\eta d)\iota_\xi = L_\eta\iota_\xi -
\iota_\eta(L_\xi - \iota_\xi d) = L_\eta\iota_\xi -\iota_\eta L_\xi + \iota_\eta\iota_\xi d =
 L_\eta\iota_\xi - L_\xi\iota_\eta + \iota_{[\xi,\eta]} + \iota_\eta\iota_\xi d$}
\begin{multline}\label{formula:twist-by-exact}
[\phi(a_1), \phi(a_2)] = [a_1+\iota_{\pi(a_1)}B, a_2+\iota_{\pi(a_2)}B]=
\\
[a_1,a_2]+L_{\pi(a_1)}\iota_{\pi(a_2)}B - L_{\pi(a_2)}\iota_{\pi(a_1)}B +
d\iota_{\pi(a_2)}\iota_{\pi(a_1)}B =
\\
[a_1,a_2] +
\iota_{[\pi(a_1),\pi(a_2)]}B +
\iota_{\pi(a_2)}\iota_{\pi(a_1)}dB = \\
\phi([a_1,a_2]) + \iota_{\pi(a_2)}\iota_{\pi(a_1)}dB
\end{multline}
shows that  $\phi = \exp(B)$ is a morphism of Leibniz algebras if and only if $B$ is closed.
\end{proof}

\begin{remark}\label{remark:iso-twist}
The calculation \eqref{formula:twist-by-exact} (combined with the equality $\pi(a_i) = \phi(a_i))$)
says that $\phi([a_1,a_2]) = [\phi(a_1),\phi(a_2)] + \iota_{\pi(\phi(a_2))}\iota_{\pi(\phi(a_1))}(-dB) =
[\phi(a_1),\phi(a_2)]_{-dB}$ (the latter operation being the $-dB$-twisted bracket on $\widehat\cA$
as defined in \ref{twist}, \eqref{formula:twisted-by-H-bracket})
shows that, for $B\in\Omega^2_X$ (not necessarily closed), the map $\exp(B)$ is a morphism
of Courant extensions $\widehat\cA_{dB}\to\widehat\cA$ (the former being the $dB$-twist of
$\widehat\cA$).
\end{remark}

\subsection{From Lie to Leibniz}\label{sssec:LtoL}
For a sheaf $\cF$ of $\cO_X$-modules let $\dual\cF$ denote
$\shHom_{\cO_X}(\cF,\cO_X)$.

Suppose that $\widehat\cA$ is a Courant extension of $\cA$.
The pairings \eqref{pairing:gA} and \eqref{pairing:gg} yield, respectively,
the maps $\mathfrak{g}(\widehat\cA) \to \dual\cA$
and $\mathfrak{g}\to\dual{\mathfrak{g}}$. Together with the projection
$\mathfrak{g}(\widehat\cA)\to\mathfrak{g}$ and the map
$\dual\cA\to\dual{\mathfrak{g}}$
adjoint to the inclusion $\mathfrak{g}\to\cA$ they fit into the diagram
\begin{equation}\label{diag:pull-back}
\begin{CD}
\mathfrak{g}(\widehat\cA) @>>> \dual\cA \\
@VVV                    @VVV  \\
\mathfrak{g} @>>> \dual{\mathfrak{g}}
\end{CD}
\end{equation}

\begin{lemma}
The diagram \eqref{diag:pull-back} is Cartesian.
\end{lemma}
\begin{proof}
The diagram \eqref{diag:pull-back} is commutative since \eqref{pairing:gg} is induced
from the restriction of \eqref{pairing:gA} to
$\mathfrak{g}(\widehat\cA)\otimes_{\cO_X}{\mathfrak{g}}$. In fact, it extends to the
morphism of short exact sequences
\[
\begin{CD}
0 @>>> \Omega^1_X @>>> \mathfrak{g}(\widehat\cA) @>>> \mathfrak{g} @>>> 0 \\
& & @VVV  @VVV  @VVV \\
0 @>>> \dual\cT_X @>>> \dual\cA @>>> \dual{\mathfrak{g}} @>>> 0
\end{CD}
\]
induced by the pairing. In particular, the map $\Omega^1_X\to\dual\cT_X$ is the
canonical isomorphism (by \eqref{ip-o}) and the claim follows.
\end{proof}

\begin{corollary}\label{cor:fib-prod}
${\mathfrak{g}}(\widehat\cA)$ is canonically isomorphic
to $\dual\cA\times_{\dual{{\mathfrak{g}}(\cA)}}{\mathfrak{g}}(\cA)$.
\end{corollary}

Since $\cA$ is transitive there is an exact
sequence
\[
0 \to {\mathfrak{g}}\to\cA\to\cT_X \to 0 \ .
\]
The map $i : \mathfrak{g}\to\cA$ and the
pairing on $\mathfrak{g}$ give rise to the maps
\[
\begin{CD}
\dual\cA @>\dual i>> \dual{\mathfrak{g}} @<{\ip}<< {\mathfrak{g}} \ .
\end{CD}
\]
Suppose in addition that $\mathfrak{g}$ is equipped with a symmetric $\cO_X$-bilinear
pairing
\[
\ip : {\mathfrak{g}}\otimes_{\cO_X}{\mathfrak{g}} \to \cO_X
\]
which is invariant under the adjoint action of $\cA$, i.e., for $a\in\cA$ and
$b,c\in{\mathfrak{g}}$
\[
\pi(a)(\langle b,c\rangle) = \langle [a,b],c\rangle + \langle b,[a,c]\rangle
\]
holds.

Let $\widehat{\mathfrak{g}} = \dual\cA\times_{\dual{\mathfrak{g}}}{\mathfrak{g}}$ and let
$\pr : \widehat{\mathfrak{g}} \to {\mathfrak{g}}$ denote the canonical projection.
A section of $\widehat{\mathfrak{g}}$ is a pair $(\dual a,b)$, where $\dual a\in\dual\cA$
and $b\in\mathfrak{g}$, which satisfies $\dual i(\dual a)(c) = \langle b,c\rangle$ for
$c\in\mathfrak{g}$.

The Lie algebra $\mathfrak{g}$ acts on $\cA$ (by the restriction of the
adjoint action) by $\cO_X$-linear endomorphisms and the map $i :
{\mathfrak{g}} \to \cA$ is a map of $\mathfrak{g}$-modules. Therefore,
$\dual\cA$ and $\dual{\mathfrak{g}}$ are $\mathfrak{g}$-modules\footnote{The
action of $b\in{\mathfrak{g}}$ on $\dual a\in\dual\cA$ is determined by
$b(\dual a)(c) = -\dual a ([b,c]) = \dual a([c,b])$ for
$c\in\cA$.} in a natural way and the map $\dual i$ is a morphism
of such. Hence, $\widehat{\mathfrak{g}}$ is a $\mathfrak{g}$-module in a
natural way and the map $\pr$ is a morphism of $\mathfrak{g}$-modules.

As a consequence, $\widehat{\mathfrak{g}}$ acquires the canonical structure of a
Leibniz algebra with the Leibniz bracket $[\widehat a,\widehat b]$ of two sections
$\widehat a,\widehat b\in\widehat{\mathfrak{g}}$ given by the formula
 $[\widehat a,\widehat b] = \pr(\widehat a)(\widehat b)$. Explicitly, for
$\dual a_1, \dual a_2\in\dual\cA$, $b_1,b_2\in\mathfrak{g}$
\[
[(\dual a_1,b_1),(\dual a_2,b_2)] = (b_1(\dual a_2),[b_1,b_2]) =
(\langle [\bullet,b_1],b_2\rangle,[b_1,b_2]) \ .
\]

We define a symmetric $\cO_X$-bilinear pairing
\[
\ip : \widehat{\mathfrak{g}}\otimes_{\cO_X}\widehat{\mathfrak{g}} \to \cO_X
\]
as the composition of $\pr\otimes\pr$ with the pairing on $\mathfrak{g}$.

The inclusion $\dual\pi:\Omega^1_X\to\dual\cA$ gives rise to the derivation
$\partial : \cO_X \to \widehat{\mathfrak{g}}$.

\begin{lemma}
The Leibniz bracket, the symmetric pairing and the derivation defined above
endow $\widehat{\mathfrak{g}}$ with the structure of a Courant extension
of $\mathfrak{g}$ (in particular, a Courant $\cO_X$-algebroid with the
trivial anchor map).

The isomorphism of Corollary \ref{cor:fib-prod} is an isomorphism
of Courant extensions of ${\mathfrak{g}}(\cA)$.
\end{lemma}
\begin{proof}
Left to the reader.
\end{proof}

Starting with a transitive Lie algebroid $\mathcal{A}$ and an
$\mathcal{A}$-invariant symmetric pairing $\ip$ on $\mathfrak{g} =
\mathfrak{g}(\mathcal{A})$ we have constructed a Courant extension
$\widehat{\mathfrak{g}}$ with the property that for any Courant
extension $\widehat{\mathcal{A}}$ which induces the pairing $\ip$ on
$\mathfrak{g}$ there is a canonical isomorphism (Corollary
\ref{cor:fib-prod})
$\widehat{\mathfrak{g}}\cong\mathfrak{g}(\widehat{\mathcal{A}})$.

In what follows we will suppress the dependence of the extension
$\widehat{\mathfrak{g}}$ on the pairing $\ip$ whenever the pairing
is fixed by the context. In this case we will use the canonical
isomorphism of Corollary \ref{cor:fib-prod} identify
$\widehat{\mathfrak{g}}$ and $\mathfrak{g}(\widehat{\mathcal{A}})$

\subsection{Leibniz extensions from connections}\label{ss:Leib-ext-conn}
Suppose that $\nabla$ is a connection on $\cA$. $\nabla$ determines
\begin{enumerate}
\item the isomorphism ${\mathfrak{g}}\bigoplus\cT_X\cong\cA$
by $(a,\xi) \mapsto i(a) + \nabla(\xi)$, where $a\in\mathfrak{g}$ and
$\xi\in\cT_X$;
\item the isomorphism $\phi_\nabla : \widehat{{\mathfrak{g}}} \to
\Omega^1_X\bigoplus{\mathfrak{g}}$ by $(\dual a,b)\mapsto
(\dual\nabla(\dual a),b)$, where $\dual a\in\dual\cA$, $b\in\mathfrak{g}$,
$\dual i(\dual a) = \langle b,\bullet\rangle$,
$\dual\nabla : \dual\cA \to \Omega^1_X$ is the transpose of $\nabla$
and $\dual i : \dual\cA \to \dual{\mathfrak{g}}$ is the transpose of $i$.
\end{enumerate}

Let $[\ ,\ ]_\nabla$ denote the Leibniz bracket on
$\Omega^1_X\bigoplus{\mathfrak{g}}$ induced by $\phi_\nabla$. Let
$c_\nabla(\ ,\ ) : {\mathfrak{g}}\otimes{\mathfrak{g}} \to
\Omega^1_X$ be (the $\Omega^1_X$-valued Leibniz cocycle) determined
by $\iota_\xi c_\nabla(a,b) = \langle [\nabla(\xi),a],b\rangle$,
where $\xi\in\cT_X$, $a,b\in\mathfrak{g}$, and the bracket is
computed in $\cA$.

\begin{lemma}\label{sssec:cocycle}

\begin{enumerate}
\item $[\ ,\ ]_\nabla$ is the extension of the Lie bracket on $\mathfrak{g}$
by the $\Omega^1_X$-valued (Leibniz) cocycle
$c_\nabla(\ ,\ )$.
\item
Suppose that $A\in\Omega^1_X\otimes_{\cO_X}{\mathfrak{g}}$. The automorphism of
$\Omega^1_X\bigoplus{\mathfrak{g}}$ defined by
\[
(\alpha, a)\mapsto(\alpha+\langle A,a\rangle,a)
\]
is the isomorphism
of Courant algebroids $(\Omega^1_X\bigoplus{\mathfrak{g}},[\ ,\ ]_\nabla) \to
(\Omega^1_X\bigoplus{\mathfrak{g}},[\ ,\ ]_{\nabla+A})$ which corresponds to the
identity map on $\widehat{\mathfrak{g}}$ under the identifications $\phi_\nabla$
and $\phi_{\nabla + A}$.
\end{enumerate}
\end{lemma}

\subsection{Courant extensions with connection}\label{ss:cext-conn}
Suppose that $\cA$ is a transitive Lie algebroid locally free of finite rank
over $\cO_X$ and $\widehat\cA$ is a Courant extension of $\cA$.

For a connection $\nabla$ on $\cA$ let $\cE_\nabla\subset\widehat\cA$ denote the ``inverse image'' of $\nabla(\cT_X)$ under the projection. Thus, $\cE_\nabla$ contains (the image of) $\Omega^1_X$ and the anchor map induces the isomorphism
$\cE_\nabla/\Omega^1_X\to\cT_X$ so that there is a short exact sequence
\begin{equation}\label{ses:E_nabla}
0\to\Omega^1_X\to\cE_\nabla\to\cT_X\to 0 \ .
\end{equation}
It is clear that the restriction of the symmetric pairing to
$\cE_\nabla$ is non-degenerate. A connection $\widehat\nabla
:\cT_X\to\widehat\cA$ on $\widehat\cA$ lifting $\nabla$ (i.e.
$\nabla$ is the composition of $\widehat\nabla$ and the projection
$\widehat\cA\to\cA$) determines a Lagrangian splitting of
\eqref{ses:E_nabla}.

Let $\cE_\nabla^\perp\subset\widehat\cA$ denote the annihilator of $\cE_\nabla$ with respect to the symmetric pairing.

\begin{lemma}\label{lemma:E_nabla_props}
\begin{enumerate}
\item $\cE_\nabla\bigcap\cE_\nabla^\perp = 0$

\item The projection $\widehat\cA\to\cA$ restricts to an isomorphism
$\cE_\nabla^\perp\to\mathfrak{g}$.

\item ${\mathfrak{g}}(\widehat\cA)$ decomposes into the orthogonal direct sum
$\Omega^1_X+\cE_\nabla^\perp$.

\item The induced isomorphism
$\Omega^1_X\oplus{\mathfrak{g}}\isomo\widehat{\mathfrak{g}}$
coincides with the one in Lemma \ref{sssec:cocycle}.

\item Suppose in addition that $\nabla$ is flat. Then, the Leibniz bracket and
the symmetric pairing restrict to a structure of a Courant algebroid
on $\cE_\nabla$.
\end{enumerate}
\end{lemma}
\begin{proof}
Since the symmetric pairing on $\cE_\nabla$  is non-degenerate
$\cE_\nabla\bigcap\cE_\nabla^\perp = 0$. Note that this means that
the natural map $\cE_\nabla\oplus\cE_\nabla^\perp\to \widehat\cA$ is
an isomorphism.

Since $\cE_\nabla^\perp
\perp\Omega^1_X(\subset\cE_\nabla)$ it follows (from \eqref{ip-o})
that the composition $\cE_\nabla^\perp\to\cA\to\cT_X$ is trivial.
Hence, the map $\cE_\nabla^\perp\to\cA$ factors through
$\cE_\nabla^\perp\to\mathfrak{g}$ which, clearly, is an isomorphism.

It follows that the composition
\[
\Omega^1_X\oplus\cE_\nabla^\perp\to
\cE_\nabla\oplus\cE_\nabla^\perp\to \widehat\cA
\]
is an isomorphism onto ${\mathfrak{g}}(\widehat\cA)$.

For $a,b\in{\mathfrak{g}}$ let $\widetilde{a}, \widetilde{b}\in\cE_\nabla^\perp$
denote their respective lifts. For $q\in\cE_\nabla$, viewed as a section
of $\widehat\cA$,
\begin{multline*}
\langle q,[\widetilde{a},\widetilde{b}]\rangle_{\widehat\cA} =
- \langle [\widetilde{a},q],\widetilde{b}\rangle_{\widehat\cA} =
- \langle\overline{[\widetilde{a},q]},b\rangle_{\mathfrak{g}} = \\
- \langle[a,\nabla(\pi(q))],b\rangle_{\mathfrak{g}} =
\langle[\nabla(\pi(q)),a],b\rangle_{\mathfrak{g}} \ ,
\end{multline*}
i.e. the bracket on $\Omega^1_X\oplus{\mathfrak{g}}$ is given precisely by the cocycle
$c_\nabla(\ ,\ )$ of Lemma \ref{sssec:cocycle}.

If $\nabla$ is flat, then $\nabla(\cT)$ is closed under the Lie bracket
in $\cA$, hence $\cE_\nabla$ is closed under the Leibniz bracket in
$\widehat\nabla$.
\end{proof}

Suppose that $\widehat\nabla$ is a connection on $\widehat\cA$ which
lifts the connection $\nabla$ on $\cA$. According to Lemma
\ref{lemma:E_nabla_props}, the curvature $c(\widehat\nabla)\in
\Omega^2_X\otimes_{\cO_X}\widehat{\mathfrak{g}}$ of the  connection
$\widehat\nabla$ (see Definition \ref{definition:curvature courant},
and Remark \ref{remark:curvature}) decomposes uniquely as
$c(\widehat\nabla) = \widetilde{c(\nabla)} +
c_{rel}(\widehat\nabla)$, where
$\widetilde{c(\nabla)}\in\Omega^2_X\otimes_{\cO_X}\cE_\nabla^\perp$
is the lift of $c(\nabla)$ (the curvature of the connection
$\nabla$) and
$c_{rel}(\widehat\nabla)\in\shHom_{\cO_X}(\bigwedge{}^2\cT_X,\Omega^1_X)=
\Omega^2_X\otimes_{\cO_X}\Omega^1_X$.

\begin{lemma}\label{lemma:rel-curv-props}
In the notations introduced above
\begin{enumerate}
\item $\iota_\xi (c_{rel}(\widehat\nabla)(\xi_1,\xi_2))=
\langle [\widehat\nabla(\xi_1),\widehat\nabla(\xi_2)],\widehat\nabla(\xi)\rangle$

\item $c_{rel}$ is totally skew-symmetric, i.e.
$c_{rel}(\widehat\nabla)\in\Omega^3_X$

\end{enumerate}
\end{lemma}
\begin{proof}
Both claims follow from the calculation
\begin{eqnarray*}
\iota_\xi (c_{rel}(\widehat\nabla)(\xi_1,\xi_2)) & = &
\langle c_{rel}(\widehat\nabla)(\xi_1,\xi_2),\widehat\nabla(\xi)\rangle \\
& = & \langle c(\widehat\nabla),\widehat\nabla(\xi)\rangle \\
& = &
\langle [\widehat\nabla(\xi_1),\widehat\nabla(\xi_2)],\widehat\nabla(\xi)\rangle \\
& = & \xi_1(\langle\widehat\nabla(\xi_2),\widehat\nabla(\xi)\rangle) -
\langle[\widehat\nabla(\xi_1),\widehat\nabla(\xi)],\widehat\nabla(\xi_2)\rangle \\
& = & - \langle[\widehat\nabla(\xi_1),\widehat\nabla(\xi)],
\widehat\nabla(\xi_2)\rangle \\
& = & - \iota_{\xi_2} (c_{rel}(\widehat\nabla)(\xi_1,\xi))
\end{eqnarray*}
(using $\langle\widetilde{c(\nabla)},\widehat\nabla(\xi)\rangle =
\langle\widehat\nabla([\xi_1,\xi_2]),\widehat\nabla(\xi)\rangle =
\langle\widehat\nabla(\xi_2),\widehat\nabla(\xi)\rangle = 0$).

\end{proof}

The choice of the lifting $\widehat\nabla$ of $\nabla$ gives rise to the
identification
\[
\Omega^1_X\oplus\mathfrak{g}\oplus\cT_X\to\widehat\cA \ ,
\]
$(\alpha,a,\xi)\mapsto i(\alpha) + \widetilde{a} + \widehat\nabla(\xi)$ which
\begin{enumerate}
\item maps the flag
$\Omega^1_X\subset\Omega^1_X\oplus\mathfrak{g}\subset\Omega^1_X\oplus
\mathfrak{g}\oplus\cT_X$ isomorphically onto the flag $\Omega^1_X\subset
\widehat{\mathfrak{g}}\subset\widehat\cA$ and

\item projects (modulo $\Omega^1_X$) to the identification
$\mathfrak{g}\oplus\cT_X\to\cA$ induced by the connection $\nabla$ on $\cA$.
\end{enumerate}

Let $\ip_{\widehat\nabla}$ and $[\ ,\ ]_{\widehat\nabla}$ denote, respectively,
the induced symmetric pairing and Leibniz bracket on
$\Omega^1_X\oplus\mathfrak{g}\oplus\cT_X$.

\begin{lemma}
$\ip_{\widehat\nabla}$ and $[\ ,\ ]_{\widehat\nabla}$ satisfy
\begin{enumerate}
\item $\langle\Omega^1_X\oplus\cT_X,\mathfrak{g}\rangle_{\widehat\nabla} = 0$;

\item the restriction of $\ip_{\widehat\nabla}$ to $\Omega^1_X\oplus\cT_X$
(respectively, $\mathfrak{g}$) is induced by the duality pairing (respectively,
is the pairing induced by the one on $\widehat\cA$);

\item the restriction of the Leibniz bracket to $\Omega^1_X\oplus\mathfrak{g}$
is given by Lemma \ref{lemma:E_nabla_props} (and Lemma \ref{sssec:cocycle});

\item $[\xi_1,\xi_2]_{\widehat\nabla} = [\nabla(\xi_1),\nabla(\xi_2)]_\cA +
c_{rel}(\xi_1,\xi_2)$

\item $[\xi,a]_{\widehat\nabla} = [\nabla(\xi),a]_{\cA} -
\langle c(\nabla)(\xi,\bullet),a\rangle_{\mathfrak{g}}$
\end{enumerate}
where $a\in\mathfrak{g}$, $\xi,\xi_i\in\cT_X$.
\end{lemma}
\begin{proof}
Only the last formula had not been proven earlier. In terms of the orthogonal
direct sum decomposition (Lemma \ref{lemma:E_nabla_props})
$\widehat\cA = \cE_\nabla^\perp + \cE_\nabla$
\[
[\widehat\nabla(\xi),\widetilde{a}] =
\widetilde{[\nabla(\xi),a]} + \Phi(\xi,a) \ ,
\]
where $\Phi(\xi,a)\in\Omega^1_X$ satisfies
\begin{eqnarray*}
\iota_\eta\Phi(\xi,a) & =& \langle\widehat\nabla(\eta),\Phi(\xi,a)\rangle \\
& = & \langle\widehat\nabla(\eta),[\widehat\nabla(\xi),\widetilde{a}]\rangle \\
& = & \xi(\langle\widehat\nabla(\eta),\widetilde{a}\rangle)-
\langle[\widehat\nabla(\xi),\widehat\nabla(\eta)],\widetilde{a}\rangle \\
& = & -\langle\widehat\nabla([\xi,\eta])+\widetilde{c(\nabla)(\xi,\eta)}+
c_{rel}(\widehat\nabla)(\xi,\eta),\widetilde{a}\rangle \\
& = & -\langle\widetilde{c(\nabla)(\xi,\eta)},\widetilde{a}\rangle \ .
\end{eqnarray*}
\end{proof}

\subsection{Construction of Courant extensions}
Suppose that $\cA$ is a transitive Lie algebroid locally free of finite rank
over $\cO_X$, $\mathfrak{g} = \mathfrak{g}(\cA)$, $\ip$ is a symmetric
$\cA$-invariant pairing on $\mathfrak{g}$.

We extend the pairing $\ip$ to a symmetric pairing on
$\Omega^1_X\oplus\mathfrak{g}\oplus\cT_X$ by the rules
\begin{enumerate}
\item $\langle\mathfrak{g},\Omega^1_X\oplus\cT_X\rangle = 0$
\item the restriction of $\ip$ to $\Omega^1_X\oplus\cT_X$ coincides with the
canonical symmetric pairing.
\end{enumerate}

For a connection $\nabla$ on $\cA$ and $H\in\Omega^3_X$ let
$\widehat\cA_{\nabla,H}$ denote $\Omega^1_X\oplus\cA$ identified with
$\Omega^1_X\oplus\mathfrak{g}\oplus\cT_X$ via $\nabla$ with the induced
symmetric pairing, denoted $\ip_{\nabla,H}$, the map (derivation)
$\partial : \cO_X\to\widehat\cA$, defined as the composition of the exterior
differentiation with the inclusion of $\Omega^1_X$, the (anchor) map
$\pi:\widehat\cA\to\cT_X$, denoting the projection onto $\cT_X$ and the operation
\[
[\ ,\ ]_{\nabla,H}:\widehat\cA_{\nabla,H}\times\widehat\cA_{\nabla,H}
\to\widehat\cA_{\nabla,H}
\]
defined by
\begin{enumerate}
\item the formulas
\begin{eqnarray*}
[\xi,df]_{\nabla,H} & = & \partial(\xi(f)) \\
\lbrack\xi,a\rbrack_{\nabla,H} & = & [\nabla(\xi),a]_\cA -
\langle c(\nabla)(\xi,\bullet),a\rangle_{\mathfrak{g}} \\
\lbrack\xi_1,\xi_2\rbrack_{\nabla,H} & = & [\xi_1,\xi_2] +
c(\nabla)(\xi_1,\xi_2) + H(\xi_1,\xi_2,\bullet) \\
\lbrack a_1, a_2\rbrack_{\nabla,H} & = & [a_1,a_2]_{\mathfrak{g}} +
\langle[\nabla(\bullet),a_1]_\cA,a_2\rangle_{\mathfrak{g}} \\
\lbrack a,\partial(f)\rbrack_{\nabla,H} & = &
\lbrack\partial(f),\partial(g)\rbrack_{\nabla,H} = 0
\end{eqnarray*}
where $f,g\in\cO_X$, $a,a_i\in\mathfrak{g}$, $\xi,\xi_i\in\cT_X$,
\item the Leibniz rule
\[
[\widehat{a},f\widehat{b}]_{\nabla,H} = [\widehat{a},\widehat{b}]_{\nabla,H}+
\pi(\widehat{a})(f)\widehat{b}
\]
\item the relation
\[
[\widehat{a},\widehat{b}]_{\nabla,H} = -[\widehat{b},\widehat{a}]_{\nabla,H}+
\partial\langle\widehat{a},\widehat{b}\rangle_{\widehat\cA_{\nabla,H}}
\]
where $f\in\cO_X$, and $\widehat{a},\widehat{b}\in\widehat\cA_{\nabla,H}$.
\end{enumerate}

\begin{remark}\label{remark:bracket-nabla-H}
It is clear from the above formulas that
\begin{enumerate}
\item $[\widehat\cA_{\nabla,H},\Omega^1_X\oplus\mathfrak{g}]_{\nabla,H}\subseteq
\Omega^1_X\oplus\mathfrak{g}$

\item the (restriction of) $[\ ,\ ]_{\nabla,H}$ to
$(\Omega^1_X\oplus\mathfrak{g})^{\otimes 2}$
coincide with the bracket $[\ ,\ ]_\nabla$ of \ref{ss:Leib-ext-conn},

\item the induced map
$\widehat\cA_{\nabla,H}\to\shEnd_\mathbb{C}(\Omega^1_X\oplus\mathfrak{g})$ is trivial on $\Omega^1_X$
hence induces the map $\cA_{\nabla,H}\to\shEnd_\mathbb{C}(\Omega^1_X\oplus\mathfrak{g})$; the latter map
is the direct sum of the Lie derivative action on $\Omega^1_X$ (through the quotient $\cT_X$) and the adjoint
action of $\cA$ on $\mathfrak{g}$.
\end{enumerate}
\end{remark}

\begin{lemma}
The Jacobiator of the operation $[\ ,\ ]_{\nabla,H}$ (see
\eqref{formula:jacobiator}), $J([\ ,\ ]_{\nabla,H}) :
\widehat\cA_{\nabla,H}^{\otimes 3} \to\widehat\cA_{\nabla,H}$
factors into the composition
\[
\widehat\cA_{\nabla,H}^{\otimes 3}
\stackrel{\pi^{\otimes 3}}{\longrightarrow}
\cT_X^{\otimes 3}
\stackrel{J_{\nabla,H}}{\longrightarrow}\Omega^1_X
\stackrel{i}{\longrightarrow}\widehat\cA_{\nabla,H} \ ,
\]
where
\[
J_{\nabla,H}(\xi_0,\xi_1,\xi_2) =
\iota_{\xi_2}\iota_{\xi_1}\iota_{\xi_0}(
-\frac12\langle c(\nabla)\wedge c(\nabla)\rangle +dH) \ .
\]
\end{lemma}
\begin{proof}
It follows from
\begin{enumerate}
\item $(\widehat\cA_{\nabla,H}\to\cA)\circ J([\ ,\ ]_{\nabla,H}) =
J([\ ,\ ]_\cA)\circ (\widehat\cA_{\nabla,H}\to\cA)^{\otimes 3}$

\item $J([\ ,\ ]_\cA) = 0$ (i.e. the Lie bracket on $\cA$ satisfies the
Jacobi identity)
\end{enumerate}
that $J([\ ,\ ]_\cA)$ takes values in $\Omega^1_X$. By Remark
\ref{remark:bracket-nabla-H} $J([\ ,\ ]_\cA)$ vanishes on
$\ker(\pi^{\otimes 3})$. This proves the first claim and, hence,
\begin{multline*}
J_{\nabla,H}(\xi_0,\xi_1,\xi_2) = J([\ ,\ ]_{\nabla,H})(\xi_0,\xi_1,\xi_2) = \\
[\xi_0,[\xi_1,\xi_2]+c(\nabla)(\xi_1,\xi_2)+\iota_{\xi_2}\iota_{\xi_1}H]_{\nabla,H}
\\
-[[\xi_0,\xi_1]+c(\nabla)(\xi_0,\xi_1)+\iota_{\xi_1}\iota_{\xi_0}H,\xi_2]_{\nabla,H}
\\
-[\xi_1,[\xi_0,\xi_2]+c(\nabla)(\xi_0,\xi_2)+\iota_{\xi_2}\iota_{\xi_0}H]_{\nabla,H}
\end{multline*}
The three summands expand, respectively, to
\begin{multline*}
[\xi_0,[\xi_1,\xi_2]+c(\nabla)(\xi_1,\xi_2)+\iota_{\xi_2}\iota_{\xi_1}H]_{\nabla,H}
= \\
[\xi_0,[\xi_1,\xi_2]]_{\nabla,H} + [\xi_0,c(\nabla)(\xi_1,\xi_2)]_{\nabla,H}
+[\xi_0,\iota_{\xi_2}\iota_{\xi_1}H]_{\nabla,H} =
\\
[\xi_0,[\xi_1,\xi_2]]+c(\nabla)(\xi_0,[\xi_1,\xi_2]) +
\iota_{[\xi_1,\xi_2]}\iota_{\xi_0}H \\
+ [\nabla(\xi_0),c(\nabla)(\xi_1,\xi_2)]_\cA -
\langle c(\nabla)(\xi_0,\bullet),c(\nabla)(\xi_1,\xi_2)\rangle
+ L_{\xi_0}\iota_{\xi_2}\iota_{\xi_1}H \ ,
\end{multline*}
\begin{multline*}
[[\xi_0,\xi_1]+c(\nabla)(\xi_0,\xi_1)+\iota_{\xi_1}\iota_{\xi_0}H,\xi_2]_{\nabla,H}
= \\
[[\xi_0,\xi_1],\xi_2]_{\nabla,H} + [c(\nabla)(\xi_0,\xi_1),\xi_2]_{\nabla,H}]
+ [\iota_{\xi_1}\iota_{\xi_0}H,\xi_2]_{\nabla,H} = \\
[[\xi_0,\xi_1],\xi_2] + c(\nabla)([\xi_0,\xi_1],\xi_2) +
\iota_{\xi_2}\iota_{[\xi_0,\xi_1]}H \\
- [\nabla(\xi_2), c(\nabla)(\xi_0,\xi_1)]_\cA +
\langle c(\nabla)(\xi_2,\bullet),c(\nabla)(\xi_0,\xi_1)\rangle
- \iota_{\xi_2}d\iota_{\xi_1}\iota_{\xi_0}H
\end{multline*}
and
\begin{multline*}
[\xi_1,[\xi_0,\xi_2]+c(\nabla)(\xi_0,\xi_2)+\iota_{\xi_2}\iota_{\xi_0}H]_{\nabla,H}
=\\
[\xi_1,[\xi_0,\xi_2]]_{\nabla,H} + [\xi_1,c(\nabla)(\xi_0,\xi_2)]_{\nabla,H}
+ [\xi_1,\iota_{\xi_2}\iota_{\xi_0}H]_{\nabla,H} = \\
[\xi_1,[\xi_0,\xi_2]] + c(\nabla)(\xi_1,[\xi_0,\xi_2]) +
\iota_{[\xi_0,\xi_2]}\iota_{\xi_1}H \\
+ [\nabla(\xi_1),c(\nabla)(\xi_0,\xi_2)]_\cA -
\langle c(\xi_1,\bullet),c(\nabla)(\xi_0,\xi_2)\rangle
+ L_{\xi_1}\iota_{\xi_2}\iota_{\xi_0}H
\end{multline*}
Summing these up one obtains
\begin{multline*}
J_{\nabla,H}(\xi_0,\xi_1,\xi_2) = \\
-\langle c(\nabla)(\xi_0,\bullet), c(\nabla)(\xi_1,\xi_2)\rangle
-\langle c(\nabla)(\xi_0,\xi_1),c(\nabla)(\xi_2,\bullet)\rangle \\
+\langle c(\nabla)(\xi_1,\bullet),c(\nabla)(\xi_0,\xi_2)\rangle
+\iota_{\xi_2}\iota_{\xi_1}\iota_{\xi_0}dH = \\
\iota_{\xi_2}\iota_{\xi_1}\iota_{\xi_0}(
-\frac12\langle c(\nabla)\wedge c(\nabla)\rangle +dH)
\end{multline*}
\end{proof}

\begin{corollary}
In the notations introduced above, $\widehat\cA_{\nabla,H}$ is a Courant
extension of $\cA$ if and only if
\[
dH = \frac12\langle c(\nabla)\wedge c(\nabla)\rangle \ .
\]
If the latter condition is fulfilled the canonical connection $\widehat\nabla$
on $\widehat\cA_{\nabla,H}$ (given the inclusion of the direct summand)
satisfies $c_{rel}(\widehat\nabla) = H$.
\end{corollary}

\begin{corollary}\label{cor:iso-conn-CExt}
In the notations of \ref{ss:cext-conn},
$dc_{rel}(\widehat\nabla)=\displaystyle\frac12\langle c(\nabla)\wedge c(\nabla)\rangle$ and the map
\begin{equation}\label{map:-nabla-hat}
\phi_{\widehat\nabla} :
\widehat\cA_{\nabla,c_{rel}(\widehat\nabla)}\to\widehat\cA
\end{equation}
induced by $\widehat\nabla$ is an isomorphism of Courant extensions of $\cA$.
\end{corollary}

\subsection{Change of connection}
Suppose that $\cA$ is a Lie algebroid, $\ip$ is a $\cA$-invariant pairing
on $\mathfrak{g} :=\mathfrak{g}(\cA)$, $\nabla$ is a connection on $\cA$,
$H\in\Omega^3_X$ satisfies $dH=\displaystyle\frac12\langle c(\nabla)\wedge c(\nabla)\rangle$.

Suppose that $\nabla^\prime$ is another connection on $\cA$. Then,
the formula
\[
\widehat\nabla^\prime(\xi) = -\frac12\langle A(\xi),A(\bullet)\rangle
+ A(\xi) +\xi \ ,
\]
where $A = \nabla^\prime - \nabla\in\Omega^1_X\otimes_{\cO_X}\mathfrak{g}$,
determines a connection on $\widehat\cA_{\nabla,H}$ which induced the connection
$\nabla^\prime$ on $\cA$.

\begin{lemma}\label{lemma:rel-curv-change-conn}
In the notations as above,
\[
c_{rel}(\widehat\nabla^\prime) = H +\langle c(\nabla)\wedge A\rangle
+ \frac12\langle[\nabla,A]\wedge A\rangle + \frac16\langle[A,A]\wedge A\rangle
\]
\end{lemma}
\begin{proof}
\begin{multline*}
[\widehat\nabla^\prime(\xi_0),\widehat\nabla^\prime(\xi_1)]_{\nabla,H} = \\
[-\frac12\langle A(\xi_0),A(\bullet)\rangle + A(\xi_0) +\xi_0,
-\frac12\langle A(\xi_1),A(\bullet)\rangle + A(\xi_1) +\xi_1
]_{\nabla,H} = \\
\frac12\iota_{\xi_1}d\langle A(\xi_0),A(\bullet)\rangle
+ [A(\xi_0),A(\xi_1)] + \langle[\nabla(\bullet),A(\xi_0)],A(\xi_1)\rangle
\\ +
[A(\xi_0),\nabla(\xi_1)] + \langle c(\nabla)(\xi_1,\bullet),A(\xi_0)\rangle
-\frac12L_{\xi_0}\langle A(\xi_1),A(\bullet)\rangle
\\
+ [\nabla(\xi_0),A(\xi_1)] -\langle c(\nabla)(\xi_0,\bullet),A(\xi_1)\rangle
+[\xi_0,\xi_1] + c(\nabla)(\xi_0,\xi_1) +\iota_{\xi_1}\iota_{\xi_0}H
\end{multline*}
Pairing the result with $\widehat\nabla^\prime(\xi_2) =
-\dfrac12\langle A(\xi_2),A(\bullet)\rangle
+ A(\xi_2) +\xi_2$ gives

\begin{multline*}
\iota_{\xi_2} c_{rel}(\widehat\nabla^\prime)(\xi_0,\xi_1) =
\langle [\widehat\nabla^\prime(\xi_0),\widehat\nabla^\prime(\xi_1)]_{\nabla,H},
\widehat\nabla^\prime(\xi_2)\rangle =
\\
+\frac12\iota_{\xi_2}\iota_{\xi_1}d\langle A(\xi_0),A(\bullet)\rangle
+ \langle[A(\xi_0),A(\xi_1)],A(\xi_2)\rangle
\\
+\langle[\nabla(\xi_2),A(\xi_0)],A(\xi_1)\rangle
+\langle[A(\xi_0),\nabla(\xi_1)],A(\xi_2)\rangle
+\langle c(\nabla)(\xi_1,\xi_2),A(\xi_0)\rangle
\\
-\frac12\iota_{\xi_2} L_{\xi_0}\langle A(\xi_1),A(\bullet)\rangle
+ \langle[\nabla(\xi_0),A(\xi_1)], A(\xi_2)\rangle
\\
-\langle c(\nabla)(\xi_0,\xi_2),A(\xi_1)\rangle
-\frac12\langle A(\xi_2),A([\xi_0,\xi_1])\rangle
\\
+\langle c(\nabla)(\xi_0,\xi_1),A(\xi_2)\rangle
+\iota_{\xi_2}\iota_{\xi_1}\iota_{\xi_0}H
\end{multline*}
The identities
\begin{multline*}
\frac12\iota_{\xi_2}\iota_{\xi_1}d\langle A(\xi_0),A(\bullet)\rangle
+\langle[\nabla(\xi_2),A(\xi_0)],A(\xi_1)\rangle
+\langle[A(\xi_0),\nabla(\xi_1)],A(\xi_2)\rangle
\\
-\frac12\iota_{\xi_2} L_{\xi_0}\langle A(\xi_1),A(\bullet)\rangle
+ \langle[\nabla(\xi_0),A(\xi_1)], A(\xi_2)\rangle
-\frac12\langle A(\xi_2),A([\xi_0,\xi_1])\rangle
\\
= \iota_{\xi_2}\iota_{\xi_1}\iota_{\xi_0}
\frac12\langle[\nabla,A]\wedge A\rangle\ ,
\end{multline*}
\begin{multline*}
-\langle c(\nabla)(\xi_2,\xi_1),A(\xi_0)\rangle
- \langle c(\nabla)(\xi_0,\xi_2),A(\xi_1)\rangle \\
+\langle c(\nabla)(\xi_0,\xi_1),A(\xi_2)\rangle =
\iota_{\xi_2}\iota_{\xi_1}\iota_{\xi_0}\langle c(\nabla)\wedge A\rangle \ ,
\end{multline*}
and
\[
\langle[A(\xi_0),A(\xi_1)],A(\xi_2)\rangle =
\frac12\langle[A,A](\xi_0,\xi_1),A(\xi_2)\rangle =
\iota_{\xi_2}\iota_{\xi_1}\iota_{\xi_0}\frac16\langle[A,A]\wedge A\rangle
\]
give
\[
c_{rel}(\widehat\nabla^\prime) =
H + (\langle c(\nabla)\wedge A\rangle
+ \frac12\langle[\nabla,A]\wedge A\rangle + \frac16\langle[A,A]\wedge A\rangle)
\ .
\]
\end{proof}

\begin{notation}
Suppose that $\cA$ is a transitive Lie algebroid and $\ip$ is an $\cA$-invariant
symmetric pairing on $\mathfrak{g}(\cA)$. For connections $\nabla,\nabla^\prime$
on $\cA$ let
\begin{equation}\label{notation:secondary-pont}
\cP(\nabla,\nabla^\prime) \stackrel{def}{=}
\langle c(\nabla)\wedge A\rangle
+\frac12\langle[\nabla,A]\wedge A\rangle + \frac16\langle[A,A],A\rangle \ ,
\end{equation}
where $A = \nabla^\prime - \nabla\in\Omega^1_X\otimes_{\cO_X}\mathfrak{g}$.
\end{notation}

\begin{lemma}\label{lemma:iso-change-conn-CExt}
The isomorphism
\[
\widehat\cA_{\nabla^\prime,H +\cP(\nabla,\nabla^\prime)}
\to
\widehat\cA_{\nabla,H}
\]
of Corollary \ref{cor:iso-conn-CExt} is given by
\begin{multline}\label{map:iso-change-conn-CExt}
\alpha + a + \xi \mapsto \alpha-\langle a,A(\bullet)\rangle+
a + \widehat\nabla^\prime(\xi) = \\
(\alpha -\langle a,A(\bullet)\rangle -\frac12\langle A(\xi),A(\bullet)\rangle)
+ (a + A(\xi)) +\xi
\end{multline}
\end{lemma}
\begin{proof}
In the orthogonal decomposition $\widehat\cA_{\nabla,H} = \cE_{\nabla^\prime}
\bigoplus\cE_{\nabla^\prime}^\perp$, the summands are given by
\[
\cE_{\nabla^\prime} = \left\lbrace
\alpha + A(\xi) + \xi\ \vert\alpha\in\Omega^1_X,\ \xi\in\cT_X
\right\rbrace
\]
and
\[
\cE_{\nabla^\prime}^\perp = \left\lbrace
-\langle a,A(\bullet)\rangle +a\ \vert a\in\mathfrak{g}\right\rbrace \ .
\]
\end{proof}
\begin{notation}
Suppose that $\cA$ is a transitive Lie algebroid and $\ip$ is an $\cA$-invariant
symmetric pairing on $\mathfrak{g}(\cA)$. For connections $\nabla,\nabla^\prime$
on $\cA$ we denote by $\phi(\nabla,\nabla^\prime)$ the
isomorphism of Lemma \ref{lemma:iso-change-conn-CExt}, where
$A=\nabla^\prime-\nabla$.
\end{notation}

\begin{lemma}\label{lemma:three-conns}
Suppose that $\nabla, \nabla^\prime, \nabla^{\prime\prime}$ are connections
on $\cA$. Then,
\begin{equation}\label{formula:triple-composition}
\phi(\nabla,\nabla^\prime)\circ\phi(\nabla^\prime,\nabla^{\prime\prime})\circ
\phi(\nabla^{\prime\prime},\nabla) = \exp(-\frac12\langle A\wedge A^\prime\rangle)\ ,
\end{equation}
where $A=\nabla^\prime-\nabla$, $A^\prime=\nabla^{\prime\prime}-\nabla^\prime$
and $\exp(\bullet)$ is as in \eqref{map-exp}.
\end{lemma}
\begin{proof}
It is clear that the left hand side of \eqref{formula:triple-composition}
is of the form $\exp(B)$ for suitable $B\in\Omega^2_X$, i.e. its value on
an element of $\widehat\cA_{\nabla,H}$ depends only on the projection of that
element to $\cT_X$. Hence, it suffices to calculate the left hand side of
\eqref{formula:triple-composition} the case $\alpha=0$, $a=0$
in the notations introduced above.

Since $\nabla - \nabla^{\prime\prime} = -(A+A^\prime)$, the formula
\eqref{map:iso-change-conn-CExt} gives
\[
\phi(\nabla^{\prime\prime},\nabla)(\xi) =
-\dfrac12\langle A(\xi)+A^\prime(\xi),A(\bullet)+A^\prime(\bullet)\rangle
- (A(\xi) + A^\prime(\xi)) + \xi
\]
\begin{multline*}
\phi(\nabla^\prime,\nabla^{\prime\prime})\circ\phi(\nabla^{\prime\prime},\nabla)=
- \dfrac12\langle A(\xi)+A^\prime(\xi),A(\bullet)+A^\prime(\bullet)\rangle  \\
+
\langle A(\xi)+ A^\prime(\xi),A^\prime(\bullet)\rangle -
\frac12\langle A^\prime(\xi),A^\prime(\bullet)\rangle
- (A(\xi)+ A^\prime(\xi))+A^\prime(\xi) + \xi = \\
-\frac12\langle A(\xi),A(\bullet)\rangle -\frac12\langle A^\prime(\xi),A(\bullet)\rangle
+\frac12\langle A(\xi),A^\prime(\bullet)\rangle - A(\xi) + \xi
\end{multline*}
and, finally,
\begin{multline*}
\phi(\nabla,\nabla^\prime)\circ\phi(\nabla^\prime,\nabla^{\prime\prime})\circ
\phi(\nabla^{\prime\prime},\nabla) = -\frac12\langle A(\xi),A(\bullet)\rangle
-\frac12\langle A^\prime(\xi),A(\bullet)\rangle \\
+\frac12\langle A(\xi),A^\prime(\bullet)\rangle
+\langle A(\xi),A(\bullet)\rangle - \frac12\langle A(\xi),A(\bullet)\rangle - A(\xi)
+A(\xi) + \xi = \\
\iota_\xi\left(-\frac12\langle A\wedge A^\prime\rangle\right) + \xi
\end{multline*}
as desired.
\end{proof}

\subsection{Exact Courant algebroids}

\begin{definition}
The Courant algebroid $\cQ$ is called {\em exact} if the anchor
map $\pi: \overline{\cQ} \to \cT_X$ is an isomorphism. Equivalently,
an exact Courant algebroid is a Courant extension of the Lie algebroid
$\cT_X$.
\end{definition}

We denote the stack of exact Courant $\cO_X$-algebroids by
$\ECA_{\cO_X}$. As was pointed out in Remark \ref{remark:CEXT-gpd},
$\ECA_{\cO_X}$ is a stack in groupoids.

For an exact Courant algebroid $\cQ$ the exact sequence \eqref{ses:assoc-Lie-alg}
takes the shape
\begin{equation}\label{ses:ECA}
0 \to \Omega^1_X \to \cQ \to \cT_X \to 0 \ .
\end{equation}
An isotropic splitting of \eqref{ses:ECA} (i.e. a connection on $\cQ$) is
necessarily Lagrangian.

\begin{lemma}
Suppose that $\cQ$ is an exact Courant algebroid.
\begin{enumerate}
\item For $\nabla\in\cC(\cQ)$ and $\omega\in\Omega^2_X$ the map $\cT_X \to \cQ$
defined by $\xi\mapsto\nabla(\xi)+i(\iota_\xi\omega)$ is a connection.
\item $\nabla\mapsto\nabla+\omega$ (where $\nabla+\omega$ is the connection
defined by the formula above) is an action of (the sheaf of groups) $\Omega^2_X$
on $\cC(\cQ)$ which endows the latter with the structure of an $\Omega^2_X$-torsor.
\end{enumerate}
\end{lemma}
\begin{proof}
The difference of two sections of the anchor map $\cQ \to \cT_X$
is a map $\cT_X \to \Omega^1_X$ or, equivalently, a section of
$\Omega^1_X\otimes_{\cO_X}\Omega^1_X$. The difference of two
isotropic sections gives rise to a skew-symmetric tensor, i.e.
a section of $\Omega^2_X$. Indeed, suppose that $\nabla$ is a connection
and $\phi : \cT_X \to \Omega^1_X$, for example, $\phi(\xi) = \iota_\xi\omega$,
where $\omega\in\Omega^2_X$. Then, for $\xi,\eta\in\cT_X$
\begin{multline*}
\langle (\nabla +\phi)(\xi),(\nabla + \phi)(\eta)\rangle = \\
\langle\nabla(\xi),\nabla(\eta)\rangle + \langle\nabla(\xi),\phi(\eta)\rangle
+ \langle\phi(\xi),\nabla(\eta)\rangle + \langle\phi(\xi),\phi(\eta)\rangle
\end{multline*}
where we use $\phi$ to denote $i\circ\phi$. Since $\nabla$ and $i$ are isotropic,
$\nabla+\phi$ is isotropic if and only if
$\langle\nabla(\xi),\nabla(\eta)\rangle + \langle\nabla(\xi),\phi(\eta)\rangle = 0$
which is equivalent to $\iota_\xi\phi(\eta) = -\iota_\eta\phi(\xi)$
by \eqref{ip-o}.

By Lemma \ref{lemma:loc-conn}, $\cC(\cQ)$ is locally non-empty, hence a torsor.
\end{proof}

For a connection $\nabla$ on an exact Courant algebroid curvature coincides with
relative curvature, is a differential $3$-form by Lemma \ref{lemma:rel-curv-props}, and will be denoted
$c(\nabla)$.

\begin{lemma}\label{lemma:curv}
\begin{enumerate}
\item The curvature form $c(\nabla)$ is closed.

\item For $\alpha\in\Omega^2_X$, $c(\nabla+\alpha) =
c(\nabla) + d\alpha$.
\end{enumerate}
\end{lemma}
\begin{proof}
The curvature form is closed since its derivative is the form $\pont_{\cT_X}$ and
the latter vanishes. The second claim follows from Lemma
\ref{lemma:rel-curv-change-conn} with $A = 0$, $H = c(\nabla)$ and $B = \alpha$.
\end{proof}

\begin{corollary}
\begin{enumerate}
\item The action of $\Omega^2_X$ on $\cC(\cQ)$ restricts to an action of
$\Omega^{2,cl}_X$ on $\cC^\flat(\cQ)$.
\item $\cC^\flat(\cQ)$, if (locally) non-empty, is an $\Omega^{2,cl}_X$-torsor.
\end{enumerate}
\end{corollary}
\begin{proof}
Suppose that $\nabla$ is a flat connection and $\alpha\in\Omega^2_X$. It follows
from Lemma \ref{lemma:curv} that $\nabla+\alpha$ is flat if and only if $\alpha$
is closed.
\end{proof}

\begin{remark}
$\cC^\flat(\cQ)$ is locally non-empty if the Poincar\'e Lemma is satisfied. This
the case in the $C^\infty$ and the analytic setting.
\end{remark}

\begin{example}\label{example:the-ECA}
The sheaf $\Omega^1_X\oplus\cT_X$ endowed with the canonical symmetric bilinear
form deduced from the duality pairing carries the canonical structure of
an exact Courant algebroid with the obvious anchor map and the derivation,
and the unique Leibniz bracket, such that the inclusion of $\cT_X$ is a flat
connection. We denote this exact Courant algebroid by $\cQ_0$.

We leave it as an exercise for the reader to write down the
explicit formula for the Leibniz bracket. The skew-symmetrization
of this bracket was discovered by T.~Courant (\cite{C}) and is usually
referred to as ``the Courant bracket''.
\end{example}

\subsection{Classification of Exact Courant algebroids}

\begin{definition}
A $\Omega^2_X\to\Omega^{3,cl}_X$-torsor is a pair $(\cC,c)$, where
$\cC$ is a $\Omega^2_X$-torsor and $c$ is a map $\cC \to
\Omega^{3,cl}_X$ which satisfies $c(s+\alpha)=c(s)+d\alpha$. A
morphism of $\Omega^2_X\to\Omega^{3,cl}_X$-torsors is a morphism of
$\Omega^2_X$-torsors which commutes with the respective maps to
$\Omega^{3,cl}_X$.
\end{definition}

Suppose that $\cQ$ is an exact Courant $\cO_X$-algebroid. The assignment
$\nabla\mapsto c(\nabla)$ gives rise to the morphism
\[
c : \cC(\cQ) \to \Omega^{3,cl}_X \ .
\]
which satisfies $c(\nabla+\alpha) = c(\nabla) + d\alpha$ by Lemma
\ref{lemma:curv}. Thus, the pair $(\cC(\cQ),c)$ is a
$(\Omega^2_X\to\Omega^{3,cl}_X)$-torsor.

\begin{lemma}\label{lemma:ECA-equiv-tors}
The correspondence $\cQ\mapsto (\cC(\cQ),c)$ establishes an
equivalence
\begin{equation}\label{map:ECA-equiv-tors}
\ECA_{\cO_X}\longrightarrow
(\Omega^2_X\to\Omega^{3,cl}_X)-\text{torsors} \ .
\end{equation}
\end{lemma}
\begin{proof}
It is clear that the association $\cQ\mapsto (\cC(\cQ),c)$ determines
a functor. We construct a quasi-inverse to the latter.

Suppose that $(\cC,c)$ is a $\Omega^2_X\to\Omega^{3,cl}_X$-torsor.
We associate to it the exact Courant algebroid which is the $(\cC,c)$-twist
of the Courant algebroid $\cQ_0$ of Example \ref{example:the-ECA} and is
constructed as follows.

The underlying extension of $\cT_X$ by $\Omega^1_X$ is the
$\cC$-twist $\cQ_0^\cC$ of the trivial extension $\cQ_0=
\Omega^1_X\oplus\cT_X$, i.e. $\cQ_0^\cC =
\mathcal{C}\times_{\Omega^2_X}\mathcal{Q}_0$. Since the action of
$\Omega^2_X$ on $\cQ_0$ preserves the symmetric pairing, it follows
that $\cQ_0^\cC$ has the induced symmetric pairing.

The Leibniz bracket on $\cQ_0^\cC$ is defined by the formula
\[
[(s_1,q_1), (s_2,q_2)] = (s_1,[q_1,q_2+\iota_{\pi(q_2)}
(s_1-s_2)]_0 + \iota_{\pi(q_1)\wedge\pi(q_2)}c(s_1)) \ ,
\]
where $s_i\in\cC$, $q_i\in\cQ_0$, $s_1-s_2\in\Omega^2_X$ is the unique
form such that $s_1=s_2 + (s_1-s_2)$ and $[\ ,\ ]_0$ denotes the (Courant)
bracket on $\cQ_0$.

Next, we verify that the bracket is, indeed, well-defined on $\cQ_0^\cC$,
i.e. is independent of the choice of particular representatives.
For $\omega\in\Omega^2_X$
\begin{multline*}
[(s_1,q_1), (s_2,q_2 + \iota_{\pi(q_2)}\omega)] = \\
(s_1,[q_1,q_2+
\iota_{\pi(q_2)}\omega + \iota_{\pi(q_2)}(s_1-s_2)]_0 +
\iota_{\pi(q_1)\wedge\pi(q_2)}c(s_1)) = \\
(s_1,[q_1,q_2+ \iota_{\pi(q_2)}(s_1-(s_2-\omega))]_0 +
\iota_{\pi(q_1)\wedge\pi(q_2)}c(s_1)) = \\
[(s_1,q_1), (s_2-\omega,q_2 )]
\end{multline*}
(using $\pi(q_2) = \pi(q_2+ \iota_{\pi(q_2)}\omega)$). This shows that the bracket
is well-defined in the second variable and we can assume that $s_2 = s_1$ after
modifying $q_2$, in which case the bracket is given by the simplified formula
\[
[(s,q_1), (s,q_2)] = (s,[q_1,q_2]_0 + \iota_{\pi(q_1)\wedge\pi(q_2)}c(s))
\]
Since $\iota_{\pi(q_1)\wedge\pi(q_2)}c(s))$ is skew-symmetric in $q_1$ and $q_2$
it follows that
\begin{eqnarray*}
[(s,q_1), (s,q_2)] + [(s,q_2), (s,q_1)] & = & (s,[q_1,q_2]_0 + [q_2,q_1]_0) \\
& = & (s,d\langle q_1,q_2\rangle) \ .
\end{eqnarray*}
This shows that the symmetrized bracket is well-defined (in both variables). Since,
as we established earlier, the bracket is well-defined in the second variable, it
follows that it is well-defined (in both variables) and, moreover, satisfies
\eqref{ip-symm}. Since $c(s)$ is a closed form, the bracket satisfies the Jacobi
identity. We leave the remaining verifications to the reader.
\end{proof}

Pairs $(\cQ,\nabla)$, where $\cQ\in\ECA_{\cO_X}$ and
$\nabla$ is a connection on $\cQ$, with morphisms of pairs
defined as morphism of algebroids which commute with respective
connections give rise to a stack which we denote $\ECA\nabla_{\cO_X}$.
It is clear that $\ECA\nabla_{\cO_X}$ is a stack in groupoids.
Note that the pair $(\cQ,\nabla)$ has no non-trivial
automorphisms.

The assignment $(\cQ,\nabla)\mapsto c(\nabla)$
gives rise to the morphism of stacks
\begin{equation}\label{map:curv}
c:\ECA\nabla_{\cO_X} \to \Omega^{3,cl}_X \ ,
\end{equation}
where $\Omega^{3,cl}_X$ is viewed as discrete, i.e. the only
morphisms are the identity maps.

\begin{lemma}\label{lemma:curv-equiv}
The morphism \eqref{map:curv} is an equivalence.
\end{lemma}
\begin{proof}
The quasi-inverse associates to $H\in\Omega^{3,cl}_X$ the
$H$-twist (see \ref{twist}) $\cQ_H$ of the algebroid $\cQ_0$ of Example
\ref{example:the-ECA}. The obvious connection on  $\cQ_H$ has curvature $H$.
\end{proof}

\section{Linear algebra}

\subsection{$\ECA$ as a vector space}
Let $\EXT^1_{\cO_X}(\cT_X,\Omega^1_X)$ denote the stack
of extensions of $\cT_X$ by $\Omega^1_X$ (in the category of
$\cO_X$-modules). The passage from an exact Courant algebroid
to the associated extension as above gives rise to the faithful
functor
\[
\ECA_{\cO_X} \to \EXT^1_{\cO_X}(\cT_X,\Omega^1_X) \ .
\]

The morphism (of complexes) $(\Omega^2_X\to\Omega^{3,cl}_X)
\to \shHom_{\cO_X}(\cT_X,\Omega^1_X)$ defined by
$\Omega^2_X\ni\omega\mapsto(\xi\mapsto\iota_\xi\omega)$ induces
the ``change of the structure group'' functor
\[
(\Omega^2_X\to\Omega^{3,cl}_X)-\text{torsors} \to
\shHom_{\cO_X}(\cT_X,\Omega^1_X)-\text{torsors}
\]
and the diagram
\[
\begin{CD}
\ECA_{\cO_X} @>>> \EXT^1_{\cO_X}(\cT_X,\Omega^1_X) \\
@V\eqref{map:ECA-equiv-tors}VV          @VVV \\
(\Omega^2_X \to \Omega^{3,cl}_X)-\text{torsors} @>>>
\shHom_{\cO_X}(\cT_X,\Omega^1_X)-\text{torsors}
\end{CD}
\]
is, clearly, commutative: both compositions (from the upper left
corner to the lower right corner) consist of taking the torsor of
(locally defined) splittings of the underlying extension of
$\mathcal{T}_X$ by $\Omega^2_X$. Note that, with the exception of
the upper left corner, all stacks in the above diagram have
canonical structures of stacks in ``$\mathbb{C}$-vector spaces in
categories'', and that the morphisms between them respect these
structures. Below we will explicate the structure a stack in
``$\mathbb{C}$-vector spaces in categories'' on $\ECA_{\cO_X}$ such
that the equivalence \ref{lemma:ECA-equiv-tors} as well as the
forgetful functor to $\EXT^1_{\cO_X}(\cT_X,\Omega^1_X)$ are
morphisms of such.

Namely, given exact Courant algebroids $\cQ_1,\ldots,\cQ_n$ and complex numbers
$\lambda_1,\ldots,\lambda_n$, the ``linear combination''
$\lambda_1\cQ_1+\dots+\lambda_n\cQ_n$ is as an exact Courant algebroid $\cQ$
together with a $\cO_X$-linear map of Leibniz algebras
\[
\cQ_1\times_{\cT_X}\dots\times_{\cT_X}\cQ_n \to \cQ
\]
(with respect to the componentwise bracket on
$\cQ_1\times_{\cT_X}\dots\times_{\cT_X}\cQ_n$)
which commutes with the respective projections to $\cT_X$ and satisfies
\[
(\partial_1(f_1),\ldots,\partial_n(f_n)) \mapsto \lambda_1\partial(f_1)
+\dots +\lambda_n\partial(f_n) \ ,
\]
where $f_i\in\cO_X$ and $\partial_i$ (respectively, $\partial$) is the derivation
$\cO_X \to \cQ_i$ (respectively, $\cO_X \to \cQ$).

\subsection{The action of $\ECA_{\cO_X}$}
As before, $\cA$ is a transitive Lie $\cO_X$-algebroid locally free of finite
rank over $\cO_X$, $\mathfrak{g}$ denotes $\mathfrak{g}(\cA)$, $\ip$
is an $\cO_X$-bilinear symmetric $\cA$-invariant pairing on $\mathfrak{g}$,
$\widehat{\mathfrak{g}}$ is the Courant extension of $\mathfrak{g}$ constructed in
\ref{sssec:LtoL}.

Let $\CEXT_{\cO_X}(\cA)_\ip$ denote the substack of Courant extensions of $\cA$
which induce the given pairing $\ip$ on $\mathfrak{g}$. Note that, if $\widehat\cA$
is in $\CEXT_{\cO_X}(\cA)_\ip$, then $\mathfrak{g}(\widehat\cA)$ is canonically
isomorphic to $\widehat{\mathfrak{g}}$.

Suppose that $\cQ$ is an exact Courant $\cO_X$ algebroid and
$\widehat\cA$ is a Courant extension of $\cA$. The ``translate
by $\cQ$ of $\widehat\cA$'' is a Courant extension $\cQ +\widehat\cA$
of $\cA$ together with a $\cO_X$-linear map of Leibniz
algebras
\[
\cQ\times_{\cT_X}\widehat\cA \to \cQ  + \widehat\cA
\]
which commutes with respective projections to $\cT_X$ and satisfies
\[
(\partial_{\cQ}(f),\partial_{\widehat\cA}(g)) \mapsto
\partial(f) + \partial(g) \ ,
\]
where $f,g\in\cO_X$, $\partial_{\widehat\cA}$, $\partial_{\cQ}$,
$\partial$ are the derivations of $\widehat\cA$, $\cQ$, and
$\cQ + \widehat\cA$ respectively. In other words,  $\cQ +\widehat\cA$
the push-out of $\cQ\times_{\cT_X}\widehat\cA$ by the
addition map $\Omega_X^1\times\Omega_X^1\xrightarrow{+}\Omega^1_X$. Thus,
a section of $\cQ + \widehat\cA$ is represented by a pair $(q,a)$
with $a\in\widehat\cA$ and $q\in\cQ$ satisfying $\pi(a) =
\pi(q)\in\cT_X$. Two pairs as above are equivalent if their
(componentwise) difference is of the form $(\alpha,-\alpha)$ for
some $\alpha\in\Omega^1_X$.

For $a_i\in\widehat\cA$, $q_i\in\cQ$ with $\pi(a_i)=\pi(q_i)$ let
\begin{equation}\label{formulas:sum}
\begin{array}{rcl}
[(q_1,a_1),(q_2,a_2)] & = & ([q_1,q_2],[a_1,a_2]),\\
\langle(q_1,a_1),(q_2,a_2)\rangle & = &
\langle q_1,q_2\rangle + \langle a_1,a_2\rangle
\end{array}
\end{equation}
These operations are easily seen to descend to $\cQ + \widehat\cA$.
The derivation $\partial : \cO_X \to \cQ + \widehat\cA$ is defined as
the composition
\begin{equation}\label{map:der-sum}
\cO_X\xrightarrow{\Delta}\cO_X\times\cO_X\xrightarrow{\partial\times\partial}
\cQ\times_{\cT_X}\widehat\cA \to \cQ + \widehat\cA \ .
\end{equation}

\begin{lemma}\label{lemma:ECA-action-trans}
\begin{enumerate}
\item
The formulas \eqref{formulas:sum} and the map \eqref{map:der-sum}
determine a structure of Courant extension of $\cA$ on $\cQ +
\widehat\cA$.

\item The map ${\mathfrak{g}}(\widehat\cA) \to \cQ +
\widehat\cA$ defined by $a\mapsto(0,a)$ induces an isomorphism
${\mathfrak{g}}(\cQ + \widehat\cA)\isomo{\mathfrak{g}}(\widehat\cA)$ of
Courant extensions of ${\mathfrak{g}}(\cA)$ (by $\Omega^1_X$).

\item
Suppose that $\widehat\cA^{(1)}$, $\widehat\cA^{(2)}$ are in
$\CEXT_{\cO_X}(\cA)_\ip$. Then, there exists a unique $\cQ$ in $\ECA_{\cO_X}$,
such that $\widehat\cA^{(2)}= \cQ + \widehat\cA^{(1)}$.
\end{enumerate}
\end{lemma}
\begin{proof}
We leave the verification of the first claim to the reader.

Let $\cQ$ denote the quotient of
$\widehat\cA^{(2)}\times_{\cA}\widehat\cA^{(1)}$ by the diagonally
embedded copy of $\widehat{{\mathfrak{g}}}$. Then, $\cQ$ is an
extension of $\cT$ by $\Omega^1_X$.

The sheaf $\widehat\cA^{(2)}\times_{\cA}\widehat\cA^{(1)}$ is a Leibniz
algebra with respect to the bracket
\[
[(q_2,q_1),(q_2',q_1')] = ([q_2,q_2'],[q_1,q_1'])
\]
and carries the symmetric pairing
\[
\langle(q_2,q_1),(q_2',q_1')\rangle = \langle q_2,q_2'\rangle -
\langle q_1,q_1'\rangle \ ,
\]
where $q_i,q_i'\in \widehat\cA^{(i)}$.

Since, for any Courant algebroid $\cP$ the adjoint action of
$\cP$ on ${\mathfrak{g}}(\cP))$ and the pairing
with the latter factor through $\overline\cP$,
it follows that the diagonally embedded copy of $\widehat{{\mathfrak{g}}}$
is a Leibniz ideal and the null-space of the pairing on
$\widehat\cA^{(2)}\times_{\cA}\widehat\cA^{(1)}$. Therefore, the
Leibniz bracket and the pairing descend to $\cQ$

The derivation
$\partial : \cO_X \to \cQ$ is given by (the image in $\cQ$ of)
\[
\partial(f) = (\partial_2(f),-\partial_1(f)) \ ,
\]
for $f\in\cO_X$, where $\partial_i$ is the derivation $\cO_X \to
\widehat\cA^{(i)}$.

The Leibniz bracket, the symmetric pairing and the derivation as
above are easily seen to define a structure of a(n exact) Courant
algebroid on $\cQ$. We claim that $\widehat\cA^{(2)}= \cQ+\widehat\cA^{(1)}$.

To this end, note that the natural embedding
$\widehat\cA^{(2)}\times_{\cA}\widehat\cA^{(1)}\hookrightarrow
\widehat\cA^{(2)}\times_{\cT_X}\widehat\cA^{(1)}$ induces the embedding
$\cQ\hookrightarrow\widehat\cA^{(2)}-\widehat\cA^{(1)}$, where the latter
is the Baer difference of extensions of $\cT_X$ by $\widehat{\mathfrak{g}}$, hence the
maps
\[
\cQ\times_{\cT_X}\widehat\cA^{(1)} \hookrightarrow
(\widehat\cA^{(2)}-\widehat\cA^{(1)})\times_{\cT_X}\widehat\cA^{(1)}
\to (\widehat\cA^{(2)}-\widehat\cA^{(1)})+\widehat\cA^{(1)}\isomo
\widehat\cA^{(2)}
\]
where the last isomorphism is the canonical isomorphism of the Baer arithmetic
of extensions of $\cT_X$ by $\widehat{\mathfrak{g}}$. We leave it to the reader to
check that the composition $\cQ\times_{\cT_X}\widehat\cA^{(1)} \to \widehat\cA^{(2)}$
induces a morphism $\cQ+\widehat\cA^{(1)}=\widehat\cA^{(2)}$ of Courant extensions
of $\cA$.
\end{proof}

\subsection{Cancellation}\label{ssec:can-triv}
Suppose that $\widehat\cA$ is a Courant
extension of $\cA$. Let $\Delta : \widehat\cA\to\widehat\cA\times_\cA\widehat\cA$
denote the diagonal embedding. By definition (see the proof of Lemma
\ref{lemma:ECA-action-trans}),
$\widehat\cA - \widehat\cA =
\widehat\cA\times_\cA\widehat\cA/\Delta(\widehat{\mathfrak g})$.
Therefore, the map $\widehat\cA\to\widehat\cA - \widehat\cA$
(the composition of the diagonal with the projection) factors through the map
$\widehat\cA/\widehat{\mathfrak g} =\cT_X\to\widehat\cA\to\widehat\cA-\widehat\cA$
which is easily seen to be a section of the projection
$\widehat\cA-\widehat\cA\to\cT_X$ and, in fact, a flat connection on
the exact Courant algebroid $\widehat\cA - \widehat\cA$. Equivalently,
there is a canonical isomorphism $\cQ_0\isomo\widehat\cA - \widehat\cA$.

More generally, suppose that $\cQ$ is an exact Courant algebroid and
$\phi : \cQ+\widehat\cA\to\widehat\cA$ is a morphism of Courant extensions
of $\cA$. The (iso)morphism $\phi$ induces induces a flat connection on $\cQ$,
i.e. an isomorphism $\cQ\isomo\cQ_0$ which is the composition
\[
\cQ_0\isomo\widehat\cA-\widehat\cA\stackrel{\phi-\id}{\longleftarrow}
\left(\cQ+\widehat\cA\right)-\widehat\cA = \cQ+\left(\widehat\cA-\widehat\cA\right)
\isomo\cQ+\cQ_0 = \cQ
\]

\subsection{The stack of Courant extensions}
Suppose $\cA$ is a transitive Lie algebroid locally free of finite rank,
${\frak g} := {\frak g}(\cA)$, $\ip$ is a $\cA$-invariant symmetric pairing
on $\frak g$. Thus, the Pontryagin class
$\pont(\cA,\ip)\in H^4(X;\tau_{\leq 4}\Omega^{\geq 2}_X)$ is defined (see
\ref{ssection:pervi-pont} of Appendix).

\begin{lemma}
The stack $\CEXT_{\cO_X}(\cA)_\ip$ is locally non-empty if and only if
locally on $X$ $\cA$ admits a connection with exact Pontryagin form.
\end{lemma}
\begin{proof}
Suppose that $\widehat\cA$ is a locally defined Courant extension of $\cA$.
By Lemma \ref{lemma:loc-conn} $\widehat\cA$ admits a (locally defined)
connection, say, $\widehat\nabla$. Let $\nabla$ denote the (locally defined)
induced connection on $\cA$. Then, according to Corollary \ref{cor:iso-conn-CExt},
$\langle c(\nabla),c(\nabla)\rangle = 2dc_{rel}(\widehat\nabla)$.

Conversely, suppose that, for a locally defined connection $\nabla$ on
$\cA$ there exists a (locally defined) form $H\in\Omega^3_X$ such that
$2dH = \langle c(\nabla),c(\nabla)\rangle$. Then, $\widehat\cA_{\nabla,H}$
 is a (locally defined) object of $\CEXT_{\cO_X}(\cA)_\ip$.
\end{proof}
\begin{corollary}
$\cA$ admits a connection with exact Pontryagin form locally on $X$
if and only if $\CEXT_{\cO_X}(\cA)_\ip$ is a $\ECA_{\cO_X}$-torsor.
\end{corollary}
\begin{proof}
Follows from Lemma \ref{lemma:ECA-action-trans}.
\end{proof}

$\ECA_{\cO_X}$-torsors are classified by
$H^2(X;\Omega^2_X\to\Omega^{3,cl}_X)$.

\begin{theorem}\label{thm:main}
Suppose that $\cA$ admits a connection with exact Pontryagin form
locally on $X$. Then, the image of the class of
$\CEXT_{\cO_X}(\cA)_\ip$ under the natural map
$H^2(X;\Omega^2_X\to\Omega^{3,cl}_X)\to H^4(X;\tau_{\leq
4}\Omega^{\geq 2}_X)$ is equal to
$-\displaystyle\frac12\pont(\cA,\ip)$, where $\pont(\cA,\ip)$ is the
first Pontryagin class of $(\cA,\ip)$.
\end{theorem}
\begin{proof}
By assumption, there exists a cover
$\cU = \left\lbrace U_i\right\rbrace_{i\in I}$ of $X$ by open sets,
connections $\nabla_i$ on $\cA_i :=\cA\vert_{U_i}$ and forms
$H_i\in\Omega^3_X(U_i)$ such that $dH_i = \cP_i :=
\displaystyle\frac12\langle c(\nabla_i)\wedge c(\nabla_i)\rangle
\in\Omega^4_X(U_i)$. Let
$A_{ij} = \nabla_j - \nabla_i\in\Omega^1_X(U_i\bigcap U_j)\otimes_{\cO_X}
{\mathfrak g}$

Let $\widehat\cA_i := \widehat\cA_{\nabla_i,H_i}$ denote the corresponding
Courant extensions with connections $\widehat\nabla_i$ as in , so that
$c_{rel}(\widehat\nabla_i) = H_i$. The Courant extensions $\widehat\cA_i$ form
a system of local trivializations of the stack $\CEXT_{\cO_X}(\cA)_\ip$.

The collection of forms
$\cP_i\in\Omega^4_X(U_i)$ (respectively, $H_i\in\Omega^3_X(U_i)$)
constitutes the cochain $\cP^{4,0}\in\check C^0(\cU;\Omega^4_X)$ (respectively,
$H\in\check C^0(\cU;\Omega^3_X)$) which satisfies $dH = \cP^{4,0}$.

Let $\cP_{ij} = \cP(\nabla_i,\nabla_j)\in\Omega^3_X(U_i\bigcap U_j)$
(defined by the formula \eqref{notation:secondary-pont}),
so that $d\cP_{ij}=
\cP_j-\cP_i$ Lemma \ref{lemma:rel-curv-change-conn} . Let
$\widehat\cP_{ij} = -H_j+H_i+\cP_{ij}$.
The forms $\widehat\cP_{ij}$ are closed:
\begin{multline*}
d\widehat\cP_{ij} = - dH_j + dH_i + d\cP(\nabla_i,\nabla_j)
= - dH_j + dH_i + \cP_j - \cP_i =  \\
-(dH_j-\cP_j) + (dH_i-\cP_i) = 0 \ .
\end{multline*}

The collection of forms $\cP_{ij}$ (respectively, $\widehat\cP_{ij}$)
constitutes the cochain $\cP^{3,1}\in\check C^1(\cU;\Omega^3_X)$ (respectively,
$\widehat\cP^{3,1}\in \check C^1(\cU;\Omega^3_X)$). These satisfy
$\widehat\cP^{3,1} = \cP^{3,1} -\check\partial H$,
$d\cP^{3,1} = \check\partial\cP^{4,0}$, $d\widehat\cP^{3,1} = 0$.

Let $\cQ_{ij}:=\cQ_{\widehat\cP_{ij}}$ be the exact Courant algebroid
with connection $\widehat\nabla_{ij}$ with $c(\widehat\nabla_{ij}) =
\widehat\cP_{ij}$. Since $H_j + \widehat\cP_{ij} = H_i + \cP_{ij}$,
there are morphisms (on $U_i\bigcap U_j$)
\[
\cQ_{ij}+\widehat\cA_j\to\widehat\cA_{\nabla_j,H_j+\cP_{ij}}\to\widehat\cA_i \ ,
\]
of which the first one is defined by the formula
$(\alpha + \xi,\beta + a + \xi)\mapsto((\alpha +\beta) + a + \xi)$
(where $\alpha,\beta\in\Omega^1_X$, $a\in{\frak g}$, $\xi\in\cT_X$), while
the second one is supplied by Corollary \ref{cor:iso-conn-CExt} (and given by
the formula of Lemma \ref{lemma:iso-change-conn-CExt}).
Let $\phi_{ij}:\cQ_{ij}+\widehat\cA_j\to\widehat\cA_i$ denote the
composition of the above maps. The composition $\phi_{ij}\circ(\id+\phi_{jk})\circ(\id+\id+\phi_{ki})$:
\begin{multline*}
(\cQ_{ij}+\cQ_{jk}+\cQ_{ki})+\widehat\cA_i =
\cQ_{ij}+(\cQ_{jk}+(\cQ_{ki}+\widehat\cA_i))\to \\
\cQ_{ij}+(\cQ_{jk}+\widehat\cA_k)\to\cQ_{ij}+\widehat\cA_j\to\widehat\cA_i
\end{multline*}
(defined on $U_i\bigcap U_j\bigcap U_k$) gives rise, according to \ref{ssec:can-triv},
to the morphism
\[
\cQ_{ij}+\cQ_{jk}+\cQ_{ki}\to\cQ_0 \ ,
\]
or, equivalently, to the flat connection $\nabla^0_{ijk}$ on $\cQ_{ij}+\cQ_{jk}+\cQ_{ki}$.

On the other hand, $\cQ_{ij}+\cQ_{jk}+\cQ_{ki}$ is canonically isomorphic to
$\cQ_{\widehat\cP_{ij}+\widehat\cP_{jk}+\widehat\cP_{ki}}$, the exact
Courant algebroid with connection $\nabla_{ijk}$ whose curvature is
$\widehat\cP_{ij}+\widehat\cP_{jk}+\widehat\cP_{ki} =
\cP_{ij}+\cP_{jk}+\cP_{ki}$.

The difference $\widehat\cP_{ijk} := \nabla_{ijk} - \nabla^0_{ijk}$ is
a $2$-form on $U_i\bigcap U_j\bigcap U_k$ which satisfies
$d\widehat\cP_{ijk} = \cP_{ij}+\cP_{jk}+\cP_{ki}$. The collection
of forms $\widehat\cP_{ijk}\in\Omega^2_X(U_i\bigcap U_j\bigcap U_k)$ forms
the cochain $\widehat\cP^{2,2}\in\check C^2(\cU;\Omega^2_X)$ which satisfies
$d\widehat\cP^{2,2} = \check\partial\widehat\cP^{3,1} =
\check\partial\cP^{3,1}$ and $\check\partial\widehat\cP^{2,2} = 0$.

After the identification of the $\cO_X$-modules underlying the Courant
extensions $\widehat\cA_i$, $\widehat\cA_j$ and $\widehat\cA_k$ with
$\Omega^1_X\oplus{\mathfrak g}\oplus\cT_X$, the composition
$\phi_{ij}\circ(\id+\phi_{jk})\circ(\id+\id+\phi_{ki})$ becomes
$\phi_{ij}\circ\phi_{jk}\circ\phi_{ki} =
\exp(-\dfrac12\langle A_{ij}\wedge A_{jk}\rangle)$, the last equality
due to Lemma \ref{lemma:three-conns}. It follows from Remark \ref{remark:iso-twist}
that $\widehat\cP_{ijk} = -\displaystyle\frac12\langle A_{ij}\wedge A_{jk}\rangle$.

Let $\cP = \cP^{4,0} - \cP^{3,1} + \widehat\cP^{2,2}$,
$\widehat\cP = -\widehat\cP^{3,1} + \widehat\cP^{2,2}$.
Let $d\pm\check\partial$ denote the differential in the total complex
$\check C^\bullet(\cU;\Omega^\bullet_X)$: for $B\in\check C^i(\cU;\Omega^j_X)$
$(d\pm\check\partial)B = dB + (-1)^j\check\partial B$. Then,
$(d\pm\check\partial)\cP = (d\pm\check\partial)\widehat\cP = 0$ and
$\widehat\cP = \cP + (d\pm\check\partial)H$, i.e. $\widehat\cP$ and $\cP$
are cohomologous cycles in $\check C^\bullet(\cU;\tau_{\leq 4}\Omega^{\geq 2}_X)$.
A comparison of $\cP$ with the first Pontryagin class $\pont(\cA,\ip)$
(calculated in \ref{ssection:pervi-pont}) reveals that
$\cP = \displaystyle\frac12\pont(\cA,\ip)$.

By construction,
the class of $\CEXT_{\cO_X}(\cA)$ is represented by $-\widehat\cP$ viewed
as a cocycle of total degree two in
$\check C^\bullet(\cU;\Omega^2_X\to\Omega^{3,cl})$ (because, by our definition,
$\cQ_{ij}$ represents $\widehat\cA_i-\widehat\cA_j$ as opposed to
$\widehat\cA_j-\widehat\cA_i$) whose image (under the ``shift by two"
isomorphism which is equal to the identity map) is $-\widehat\cP$ viewed
as a cocycle of total degree four in
$\check C^\bullet(\cU;\tau_{\leq 4}\Omega^{\geq 2}_X)$. Now $-\widehat\cP$
is cohomologous to $-\cP = -\displaystyle\frac12\pont(\cA,\ip)$ which finishes
the proof.
\end{proof}

In view of \eqref{formula:pont-ch}, Theorem \ref{thm:main} can be restated in the following way for $GL_n$-torsors
(equivalently, vector bundles). In this setting we will write $\CEXT_{\cO_X}(\cA_{\cE})_{\Tr}$ for
$\CEXT_{\cO_X}(\cA_{\cE})_\ip$ as a reminder of the origins of the canonical pairing on the Atiyah algebra of a vector
bundle.

\begin{corollary}\label{cor:main-for-vb}
Suppose that $\cE$ is a vector bundle on $X$. The class of the $\ECA_{\cO_X}$-torsor $\CEXT_{\cO_X}(\cA_{\cE})_{\Tr}$
is equal to $-\ch_2(\cE)$.
\end{corollary}

\section{Vertex algebroids}\label{section:VA}

\subsection{Vertex operator algebras}
Throughout this section we follow the notations of \cite{GMS}. The following
definitions are lifted from loc. cit.

\begin{definition}
A {\it $\mathbb{Z}_{\geq 0}$-graded vertex algebra} is a
$\mathbb{Z}_{\geq 0}$-graded $k$-module $V=\oplus\ V_i$, equipped with
a distinguished vector $\vac\in V_0$ ({\it vacuum vector})
and a family of bilinear operations
\[
_{(n)}:\ V\times V\to V,\ (a, b)\mapsto a_{(n)}b
\]
of degree $-n-1$, $n\in \mathbb{Z}$, such that
\[
\vac_{(n)}a=\delta_{n,-1}a,\ \ a_{(-1)}\vac=a,\ \ a_{(n)}\vac=0\
\text{\ if\ }n\geq 0,
\]
and
\begin{multline*}
\sum_{j=0}^\infty\ \binom{m}{j}(a_{(n+j)}b)_{(m+l-j)}c= \\
\sum_{j=0}^\infty\ (-1)^j\binom{n}{j}
\bigl\{a_{(m+n-j)}b_{(l+j)}c-(-1)^nb_{(n+l-j)}a_{(m+j)}c\bigr\}
\end{multline*}
for all $a,b,c\in V,\ m,n,l\in \mathbb{Z}$.

A morphism of vertex algebras is a map of graded $k$-modules (of degree zero)
which maps the vacuum vector to the vacuum vector and commutes with all of
the operations.
\end{definition}

Let $\cV ert$ denote the category of vertex algebras.

Let
\[
\partial^{(j)}a:=a_{(-1-j)}\vac,\ \ j\in\mathbb{Z}_{\geq 0} \ .
\]
Then, $\partial^{(j)}$ is an endomorphisms of $V$ of degree $j$ which satisfies
(see \cite{GMS})
\begin{itemize}
\item $\partial^{(j)}\vac=\delta_{j,0}\vac$,
\item $\partial^{(0)}=Id$,
\item $\partial^{(i)}\cdot\partial^{(j)}=\binom{i+j}{i}\partial^{(i+j)}$,
\item $(\partial^{(j)}a)_{(n)}b=(-1)^j\binom{n}{j}a_{(n-j)}b$, and
\item $\partial^{(j)}(a_{(n)}b)=\sum_{p=0}^j\ (\partial^{(p)}a)_{(n)}\partial^{(j-p)}b$
\end{itemize}
for all $n\in\mathbb{Z}$.

The subject of the definition below is the restriction of the structure of a vertex
algebra to the graded components of degrees zero and one.

\begin{definition}
A {\it $1$-truncated vertex algebra} is a
septuple $v=(V_0,V_1,\vac,\partial,_{(-1)},_{(0)},_{(1)})$ where
\begin{itemize}
\item $V_0, V_1$ are $k$-modules,
\item $\vac$ an element of $V_0$ ({\it vacuum vector}),
\item $\partial:\ V_0\to V_1$ a $k$-linear map,
\item $_{(i)}:\ (V_0\oplus V_1)\times (V_0\oplus V_1)\to V_0\oplus V_1$
(where $(i=-1,0,1)$) are $k$-bilinear operations of degree $-i-1$.
\end{itemize}

Elements of $V_0$ (resp., $V_1$) will be denoted $a,b,c$ (resp., $x,y,z$).
There are seven operations: $a_{(-1)}b, a_{(-1)}x, x_{(-1)}a, a_{(0)}x$,
$x_{(0)}a, x_{(0)}y$ and $x_{(1)}y$. These are required to
satisfy the following axioms:
\begin{itemize}
\item
(Vacuum)
\[
a_{(-1)}\vac=a;\ x_{(-1)}\vac=x;\ x_{(0)}\vac=0
\]
\item
(Derivation)
\begin{itemize}
\item[$Deriv_1\ \ $]
$(\partial a)_{(0)}b=0;\ (\partial a)_{(0)}x=0;\
(\partial a)_{(1)}x=-a_{(0)}x$
\item[$Deriv_2\ \ $]
$\partial(a_{(-1)}b)=(\partial a)_{(-1)}b+a_{(-1)}\partial b;\
\partial(x_{(0)}a)=x_{(0)}\partial a$
\end{itemize}
\item
(Commutativity)
\begin{itemize}
\item[$Comm_{-1}$]
$a_{(-1)}b=b_{(-1)}a;\ a_{(-1)}x=x_{(-1)}a-\partial(x_{(0)}a)$
\item[$Comm_0\ \ $]
$x_{(0)}a=-a_{(0)}x;\ x_{(0)}y=-y_{(0)}x+\partial(y_{(1)}x)$
\item[$Comm_1\ \ $]
$x_{(1)}y=y_{(1)}x$
\end{itemize}
\item
(Associativity)
\begin{itemize}
\item[$Assoc_{-1}$] $(a_{(-1)}b)_{(-1)}c=a_{(-1)}b_{(-1)}c$
\item[$Assoc_0\ \ $]
$\alpha_{(0)}\beta_{(i)}\gamma=
(\alpha_{(0)}\beta)_{(i)}\gamma+\beta_{(i)}\alpha_{(0)}\gamma,\ (\alpha,
\beta, \gamma\in V_0\oplus V_1)$
whenever the both sides are defined, i.e. the operation $_{(0)}$ is a derivation
of all of the operations $_{(i)}$.
\item[$Assoc_1\ \ $]
$(a_{(-1)}x)_{(0)}b=a_{(-1)}x_{(0)}b$
\item[$Assoc_2\ \ $]
$(a_{(-1)}b)_{(-1)}x=a_{(-1)}b_{(-1)}x+(\partial a)_{(-1)}b_{(0)}x+
(\partial b)_{(-1)}a_{(0)}x$
\item[$Assoc_3\ \ $]
$(a_{(-1)}x)_{(1)}y=a_{(-1)}x_{(1)}y-x_{(0)}y_{(0)}a$
\end{itemize}
\end{itemize}

A {\it morphism} between two $1$-truncated vertex algebras
$f:\ v=(V_0,V_1,\ldots)\to v'=(V'_0,V'_1,\ldots)$ is a pair of maps of $k$-modules
$f=(f_0,f_1),\ f_i:\ V_i\to V'_i$ such that $f_0(\vac)=\vac',\ f_1(\partial a)=
\partial f_0(a)$ and $f(\alpha_{(i)}\beta)=f(\alpha)_{(i)}f(\beta)$,
whenever both sides are defined.
\end{definition}

Let $\cV ert_{\leq 1}$ denote the category of $1$-truncated vertex algebras.
We have an obvious truncation functor
$$
t:\ \cV ert\to \cV ert_{\leq 1}
\eqno{(3.1.1)}
$$
which assigns to a vertex algebra $V=\oplus_i V_i$ the truncated algebra
$tV:=(V_0,V_1,\partial^{(1)}, _{(-1)},_{(0)},_{(1)})$.

\begin{remark}\label{rem:comm-alg}
It follows easily that the operation $_{(-1)}\ :V_0\times V_0\to V_0$
endows $V_0$ with a structure of a commutative $k$-algebra.
\end{remark}

\subsection{Vertex algebroids}
Suppose that $X$ is smooth variety over $\mathbb{C}$ (a complex manifold,
a $C^\infty$-manifold). In either case we will denoted by $\cO_X$
(respectively, $\cT_X$, $\Omega^i_X$) the corresponding structure sheaf
(respectively, the sheaf of vector fields, the sheaf of differential $i$-forms).

A vertex $\cO_X$-algebroid, as defined in this section, is, essentially,
a sheaf of $1$-truncated vertex algebras, whose degree zero component (which is
a sheaf of algebras by \ref{rem:comm-alg}) is identified with $\cO_X$.

\begin{definition}\label{definition:va}
A {\em vertex $\cO_X$-algebroid} is a sheaf of $\mathbb{C}$-vector
spaces $\cV$ with a pairing
\begin{eqnarray*}
\cO_X\otimes_\mathbb{C}\cV & \to & \cV \\
f\otimes v & \mapsto & f*v
\end{eqnarray*}
such that $1* v = v$ (i.e. a ``non-associative unital
$\cO_X$-module'') equipped with
\begin{enumerate}
\item
a structure of a Leibniz $\mathbb{C}$-algebra $[\ ,\ ] :
\cV\otimes_\mathbb{C}\cV\to \cV$

\item
a $\mathbb{C}$-linear map of Leibniz algebras $\pi : \cV\to \cT_X$
(the {\em anchor})
\item
a symmetric $\mathbb{C}$-bilinear pairing $\langle\ ,\ \rangle :
\cV\otimes_\mathbb{C}\cV\to \cO_X$
\item
a $\mathbb{C}$-linear map $\partial : \cO_X\to \cV$ such that
$\pi\circ\partial = 0$
\end{enumerate}
which satisfy
\begin{eqnarray}
f*(g*v) - (fg)*v & = & \pi(v)(f)*\partial(g) +
\pi(v)(g)*\partial(f)\label{assoc} \\
\left[v_1,f*v_2\right] & = & \pi(v_1)(f)*v_2 + f*[v_1,v_2] \label{leib}
\\
\left[v_1,v_2\right] + [v_2,v_1] & = & \partial(\langle v_1,v_2\rangle)
\label{symm-bracket}\\
\pi(f*v) & = & f\pi(v) \label{anchor-lin} \\
\langle f*v_1, v_2\rangle & = & f\langle v_1,v_2\rangle -
\pi(v_1)(\pi(v_2)(f)) \label{pairing}\\
\pi(v)(\langle v_1, v_2\rangle) & = & \langle[v,v_1],v_2\rangle +
\langle v_1,[v,v_2]\rangle \label{pairing-inv} \\
\partial(fg) & = & f*\partial(g) + g*\partial(f) \label{deriv} \\
\left[v,\partial(f)\right] & = & \partial(\pi(v)(f)) \label{bracket-o}\\
\langle v,\partial(f)\rangle & = & \pi(v)(f)\label{pairing-o}
\end{eqnarray}
for $v,v_1,v_2\in\cV$, $f,g\in\cO_X$.

A morphism of vertex $\cO_X$-algebroids is a $\mathbb{C}$-linear map of
sheaves which preserves all of the structures.
\end{definition}

We denote the category of vertex $\cO_X$-algebroids by
$\VA_{\cO_X}(X)$. It is clear that the notion vertex $\cO_X$-algebroid
is local, i.e. vertex $\cO_X$-algebroids form a stack which we denote
by $\VA_{\cO_X}$.

\subsection{From vertex algebroids to 1-truncated vertex algebras}\label{ssection:VA-to-va}
Suppose that $\cV$ is a vertex $\cO_X$-algebroid. For $f,g\in\cO_X$,
$v,w\in\cV$ let
\begin{equation}\label{op:-1}
f_{(-1)}g = fg,\ \ f_{(-1)}v = f*v,\ \ v_{(-1)}f = f*v - \partial\pi(v)(f),
\end{equation}
\begin{equation}\label{op:0}
v_{(0)}f = - f_{(0)}v = \pi(v)(f),\ \ v_{(0)}w = [v,w],
\end{equation}
\begin{equation}\label{op:1}
v_{(1)}w = \langle v,w\rangle\ .
\end{equation}

\begin{lemma}
The septuple $(\cO_X,\cV,1,\partial, _{(-1)}, _{(0)}, _{(1)})$ is a
sheaf of 1-truncated vertex operator algebras.

Conversely, if the septuple $(\cO_X,\cV,1,\partial,_{(-1)},_{(0)},_{(1)})$
is a sheaf of 1-truncated vertex operator algebras, then the formulas
\eqref{op:-1}, \eqref{op:0}, \eqref{op:1} define a structure of a vertex
$\cO_X$-algebroid on $\cV$.
\end{lemma}

\subsection{Quantization of Courant algebroids}\label{ssection:quantCa}
\begin{definition}
We call a vertex algebroid {\em commutative} if the anchor map, the Leibniz
bracket and the symmetric pairing are trivial.
\end{definition}

\begin{remark}
Commutativity implies that the $*$-operation is associative, i.e.
a commutative vertex algebroid is simply a $\cO$-module $\cE$ together
with a derivation $\partial : \cO_X\to\cE$.
\end{remark}

Suppose that $\cV$ is a family of vertex $\cO_X$-algebroids flat over
$\mathbb{C}[[t]]$, such that the vertex algebroid $\cV_0 = \cV/t\cV$ is
commutative. Let $\cV\to\cV_0 : v\mapsto\overline{v}$ denote the
``reduction modulo $t$'' map. For $f\in\cO_X$, $v\in\cV$ let
$f\overline{v} = \overline{f*v}$. This operation endows $\cV_0$ with
a structure of a module over $\cO_X$.

Since $\cV_0$ is commutative, the Leibniz bracket and the symmetric
pairing on $\cV$ and the anchor map take values in $t\cV$.
For $f\in\cO_X$, $v,v_1,v_2\in\cV$ let
\begin{equation}\label{formulas:limit}
[\overline{v_1},\overline{v_2}]_0 = \overline{\frac1t[v_1,v_2]},\ \
\langle\overline{v_1},\overline{v_2}\rangle_0 =
\overline{\frac1t\langle v_1,v_2\rangle},\ \ \pi_0(\overline{v}) =
\overline{\frac1t\pi(v)},\ \ \partial_0(f) =
\overline{\partial(f)} \ .
\end{equation}

\begin{lemma}
The formulas \eqref{formulas:limit} endow the $\cO_X$-module $\cV_0$
with the structure of a Courant algebroid with derivation
$\partial_0$, anchor map $\pi_0$, Leibniz bracket $[\ ,\ ]_0$
and symmetric pairing $\langle\ ,\ \rangle_0$.
\end{lemma}

\subsection{The associated Lie algebroid}
Suppose that $\cV$ is a vertex $\cO_X$-algebroid. Let
\begin{eqnarray*}
\Omega_\cV & \stackrel{def}{=} & \cO_X*\partial(\cO_X)\subset\cV \ ,\\
\overline\cV & \stackrel{def}{=} & \cV/\Omega_\cV \ .
\end{eqnarray*}
Note that the symmetrization of the Leibniz bracket takes values
in $\Omega_\cV$.

For $f,g,h\in\cO_X$
\[
f*(g*\partial(h)) - (fg)*\partial(h) =
\pi(\partial(h))(f)*\partial(g) + \pi(\partial(h))(g)*\partial(f)
= 0\ ,
\]
because $\pi\circ\partial = 0$. Therefore, $\cO_X*\Omega_\cV =
\Omega_\cV$, and $\Omega_\cV$ is an $\cO_X$-module. The map
$\partial : \cO_X\to\Omega_\cV$ is a derivation, hence induces
the $\cO_X$-linear map $\Omega^1_X\to\Omega_\cV$.

Since the associator of the $\cO_X$-action on $\cV$ takes values
in $\Omega_\cV$, $\overline\cV$ is an $\cO_X$-module.

For $f,g,h\in\cO_X$
\[
\pi(f\partial(g))(h) = f\pi(\partial(g))(h) = 0 \ .
\]
Therefore, $\pi$ vanishes on $\Omega_\cV$, hence, factors through the
map
\begin{equation}\label{VA-anchor}
\pi : \overline\cV\to\cT_X
\end{equation}
of $\cO_X$-modules.

For $v\in\cV$, $f,g\in\cO_X$
\begin{eqnarray*}
[v,f\partial(g)] & = & \pi(v)(f)\partial(g)+f[v,\partial(g)] \\
& = & \pi(v)(f)\partial(g)+f\partial(\pi(v)(g)) \ .
\end{eqnarray*}
Therefore, $[\cV,\Omega_\cV]\subseteq\Omega_\cV$ and the Leibniz bracket on
$\cV$ descends to the operation
\begin{equation}\label{VA-bracket}
[\ ,\ ]:\overline\cV\otimes_\mathbb{C}\overline\cV\to\overline\cV
\end{equation}
which is skew-symmetric because the symmetrization of the Leibniz
bracket on $\cV$ takes values in $\Omega_\cV$ and satisfies the Jacobi
identity because the Leibniz bracket on $\cV$ does.

\begin{lemma}
The $\cO_X$-module $\overline\cV$ with the bracket
\eqref{VA-bracket} and the anchor \eqref{VA-anchor} is a Lie
$\cO_X$-algebroid.
\end{lemma}

\subsection{Transitive vertex algebroids}
\begin{definition}
A vertex $\cO_X$-algebroid is called {\em transitive} if the anchor
map is surjective.
\end{definition}

\begin{remark}
The vertex $\cO_X$-algebroid $\cV$ is called transitive if and only if
the Lie $\cO_X$-algebroid $\overline\cV$ is.
\end{remark}

Suppose that $\cV$ is a transitive vertex $\cO_X$-algebroid. The
derivation $\partial$ induces the map
\[
i : \Omega^1_X\to\cV \ .
\]
For $v\in\cV$, $f,g\in\cO_X$
\begin{eqnarray*}
\langle v,f\partial(g)\rangle & = & f\langle v,\partial(g)\rangle
- \pi(\partial(g))\pi(v)(f) \\
& = & f\pi(v)(g) \\
& = & \iota_{\pi(v)}fdg \ .
\end{eqnarray*}
If follows that the map $i$ is adjoint to the anchor map $\pi$.
The surjectivity of the latter implies the latter implies the
injectivity of the former. Since, in addition, $\pi\circ i = 0$
the sequence
\[
0\to\Omega^1_X\xrightarrow{i}\cV\xrightarrow{\pi}\overline\cV\to 0
\]
is exact and $i$ is isotropic.

\subsection{Exact vertex algebroids}
\begin{definition}
A vertex algebroid $\cV$ is called {\em exact} if the map
$\overline\cV\to\cT_X$ is an isomorphism.
\end{definition}

\begin{notation}
We denote the stack of exact vertex $\cO_X$-algebroids by
$\EVA_X$.
\end{notation}

A morphism of exact vertex algebroids induces a morphism of
respective extensions of $\cT_X$ by $\Omega^1_X$, hence is an
isomorphism of sheaves of $\mathbb{C}$-vector spaces. It is clear that
the inverse isomorphism is a morphism of vertex
$\cO_X$-algebroids. Hence, $\EVA_X$ is a stack in groupoids.

\begin{example}\label{example:loc-pic}
Suppose that $\cT_X$ is freely generated as an $\cO_X$-module by
a locally constant subsheaf of Lie $\mathbb{C}$-subalgebras
$\tau\subset\cT_X$, i.e. the canonical map $\cO_X\otimes_\mathbb{C}\tau
\to\cT_X$ is an isomorphism.

There is a unique structure of an exact vertex $\cO_X$-algebroid on
$\cV=\Omega^1_X\bigoplus\cT_X$ such that
\begin{itemize}
\item $f*(1\otimes t) = f\otimes t$ for $f\in\cO_X$, $t\in\tau$,

\item the anchor map is given by the projection $\cV\to\cT_X$,

\item the map $\tau\to\cV$ is a morphism of Leibniz algebras,

\item the pairing on $\cV$ restricts to the trivial pairing on
$\tau$,

\item the derivation $\partial : \cO_X\to\cV$ is given by the
composition $\cO_X\xrightarrow{d}\Omega^1_X\to\cV$.
\end{itemize}

Indeed, the action of $\cO_X$ is completely determined by \eqref{assoc}:
for $f,g\in\cO_X$, $t\in\tau$,
\[
f*(g\otimes t) = f*(g*(1\otimes t)) = fg\otimes t + t(f)dg + t(g)df \ .
\]
In a similar fashion the bracket is completely determined by \eqref{leib}
and \eqref{symm-bracket}, and the pairing is determined by \eqref{pairing}.

We leave the verification of the identities \eqref{assoc} - \eqref{pairing-o}
to the reader.
\end{example}

\subsection{Vertex extensions of Lie algebroids}
Suppose that $\cA$ is a Lie $\cO_X$-algebroid.

\begin{definition}
A vertex extension of $\cA$ is a vertex algebroid $\widehat\cA$ together with
an isomorphism $\overline{\widehat\cA}=\cA$ of Lie $\cO_X$-algebroids.

A morphism of vertex extensions of $\cA$ is a morphism of vertex algebroids which
is compatible with the identifications.
\end{definition}

\begin{notation}
We denote by $\VEXT_{\cO_X}(\cA)$ the stack of vertex extensions of
$\cA$.
\end{notation}

\subsection{Vertex extensions of transitive Lie algebroids}
From now on we suppose that $\cA$ is a transitive Lie
$\cO_X$-algebroid locally free of finite rank over $\cO_X$. Let
$\mathfrak{g}=\mathfrak{g}(\cA)$

Suppose that $\widehat\cA$ is a vertex extension of $\cA$. Then, the
derivation $\partial : \cO_X\to\widehat\cA$ induces an isomorphism
$\Omega^1_X\cong \Omega_{\widehat\cA}$. The the resulting exact (by
the same argument as in Lemma \ref{lemma:exactness of the exact})
sequence
\[
0\to\Omega^1_X\to\widehat\cA\to\cA\to 0
\]
is canonically associated to the vertex extension $\widehat\cA$ of $\cA$. Since
a morphism of vertex extensions of $\cA$ induces a morphism of associated extensions
of $\cA$ by $\Omega^1_X$ it is an isomorphism of the underlying sheaves. It is clear
that the inverse isomorphism is a morphism of vertex extensions of $\cA$.

Therefore, $\VEXT_{\cO_X}(\cA)$ is a stack in groupoids.

\begin{remark}
$\VEXT_{\cO_X}(\cT_X)$ is none other than $\EVA_{\cO_X}$.
\end{remark}

Suppose that $\widehat\cA$ is a vertex extension of $\cA$. Let
$\widehat{\mathfrak{g}} = \mathfrak{g}(\widehat\cA)$ denote the
kernel of the anchor map (of $\widehat\cA$). Thus, $\widehat{\frak
g}$ is a vertex (equivalently, Courant) extension of $\mathfrak{g}$.

Analysis similar to that of \ref{subsection:CEXT} shows that
\begin{itemize}
\item the symmetric pairing on $\widehat\cA$ induces a symmetric
$\cO_X$-bilinear pairing on $\mathfrak{g}$ which is $\cA$-invariant;

\item the vertex extension $\widehat{\mathfrak{g}}$ is obtained from the
Lie algebroid $\cA$ and the symmetric $\cA$-invariant pairing on
$\mathfrak{g}$ as in \ref{sssec:LtoL}.
\end{itemize}

\subsection{The action of $\ECA_{\cO_X}$}
As before, $\cA$ is a transitive Lie $\cO_X$-algebroid locally
free of finite rank over $\cO_X$, ${\mathfrak{g}}$ denotes ${\frak
g}(\cA)$, $\ip$ is an $\cO_X$-bilinear symmetric $\cA$-invariant
pairing on $\mathfrak{g}$, $\widehat{{\mathfrak{g}}}$ is the Courant
extension of $\mathfrak{g}$ constructed in \ref{sssec:LtoL}.

\begin{notation}
Let $\VEXT_{\cO_X}(\cA)_\ip$ denote the stack of Courant
extensions of $\cA$ which induce the given pairing $\ip$ on $\mathfrak{g}$.
\end{notation}

\begin{remark}
Clearly, $\VEXT_{\cO_X}(\cA)_\ip$ is a stack in groupoids.
Note that, if $\widehat\cA$ is in $\VEXT_{\cO_X}(\cA)_\ip$, then
${\mathfrak{g}}(\widehat\cA)$ is canonically isomorphic to
$\widehat{\mathfrak{g}}$.
\end{remark}

Suppose that $\cQ$ is an exact Courant $\cO_X$ algebroid and
$\widehat\cA$ is a vertex extension of $\cA$. Let $\widehat\cA +
\cQ$ denote the push-out of $\widehat\cA\times_{\cT_X}\cQ$ by the
addition map $\Omega_X^1\times\Omega_X^1\xrightarrow{+}\Omega^1_X$. Thus,
a section of $\widehat\cA + \cQ$ is represented by a pair $(a,q)$
with $a\in\widehat\cA$ and $q\in\cQ$ satisfying $\pi(a) =
\pi(q)\in\cT_X$. Two pairs as above are equivalent if their
(componentwise) difference is of the form $(i(\alpha),-i(\alpha))$ for
some $\alpha\in\Omega^1_X$.

For $a\in\widehat\cA$, $q\in\cQ$ with $\pi(a)=\pi(q)$, $f\in\cO_X$
let
\begin{equation}\label{formulas:vert-mult}
f*(a,q) = (f*a, fq),\ \ \ \partial(f) = \partial_{\widehat\cA}(f)
+ \partial_{\cQ}(f)\ .
\end{equation}

For $a_i\in\widehat\cA$, $q_i\in\cQ$ with $\pi(a_i)=\pi(q_i)$ let
\begin{equation}\label{formulas:vert-sum}
[(a_1,q_1),(a_2,q_2)] = ([a_1,a_2],[q_1,q_2]),\ \ \
\langle(a_1,q_1),(a_2,q_2)\rangle = \langle a_1,a_2\rangle +
\langle q_1,q_2\rangle
\end{equation}
These operations are easily seen to descend to $\widehat\cA +
\cQ$.

The two maps $\Omega^1_X\to\widehat\cA + \cQ$ given by
$\alpha\mapsto (i(\alpha),0)$ and $\alpha\mapsto (0,i(\alpha))$
coincide; we denote their common value by
\begin{equation}\label{map:der-sum-vert}
i : \Omega^1_X\to\widehat\cA + \cQ \ .
\end{equation}

\begin{lemma}
The formulas \eqref{formulas:vert-mult}, \eqref{formulas:vert-sum}
and the map \eqref{map:der-sum-vert} determine a structure of vertex
extension of $\cA$ on $\widehat\cA + \cQ$. Moreover, the map
the map
${\mathfrak{g}}(\widehat\cA)\to\widehat\cA+\cQ$ defined by $a\mapsto(a,0)$
induces an isomorphism ${\mathfrak{g}}(\widehat\cA+\cQ)\isomo{\mathfrak{g}}(\widehat\cA)$
of vertex (equivalently, Courant) extensions of ${\mathfrak{g}}(\cA)$ (by $\Omega^1_X$).
\end{lemma}

\begin{lemma}\label{lemma:ECA-action-trans-vert}
Suppose that $\widehat\cA^{(1)}$, $\widehat\cA^{(2)}$ are in
$\VEXT_{\cO_X}(\cA)_\ip$. Then, there exists a unique $\cQ$ in
$\ECA_{\cO_X}$, such that $\widehat\cA^{(2)}= \widehat\cA^{(1)} +
\cQ$.
\end{lemma}
\begin{proof}
Let $\cQ$ denote the quotient of
$\widehat\cA^{(2)}\times_{\cA}\widehat\cA^{(1)}$ by the diagonally
embedded copy of $\widehat{{\mathfrak{g}}}$. Then, $\cQ$ is an
extension of $\cT$ by $\Omega^1_X$. There is a unique structure of
an exact Courant algebroid on $\cQ$ defined as in Lemma {lemma:ECA-action-trans},
such that $\widehat\cA^{(2)}=\widehat\cA^{(1)} + \cQ$.
\end{proof}

\begin{proposition}
$\EVA_X$ is a torsor under $\ECA_X$.
\end{proposition}

\begin{proof}
In view of Lemma \ref{lemma:ECA-action-trans-vert} and the equality
$\EVA_X = \VEXT_{\cO_X}(\cT_X)$ it remains to show that
$\EVA_X$ is locally non-empty.

In the analytic or $C^\infty$ case example \ref{example:loc-pic}
provides a locally defined EVA. Indeed, locally on X there exists
an (abelian) Lie $\mathbb{C}$-subalgebra $\tau$
such that $\cT_X\cong\cO_X\otimes_\mathbb{C}\tau$.

In the algebraic setting the same example shows that $\EVA_X$ is
non-empty locally in \'etale topology. Since
$H^2_{\text{\'et}}(X,\Omega^2_X\to\Omega^{3,cl}_X)$ is canonically isomorphic
to $H^2(X,\Omega^2_X\to\Omega^{3,cl}_X)$ it follows that $\EVA_X$
is non-empty Zariski-locally.
\end{proof}

\subsection{Comparison of $\ECA_{\cO_X}$-torsors}\label{ss:comp-ECA-tors}
Suppose that $\widehat\cA$ is a vertex extension of the Lie
algebroid $\cA$. Let $\ip$ denote the induced symmetric pairing
on ${\mathfrak{g}}(\cA)$.

Suppose that $\cV$ is an exact vertex algebroid. Let $\widehat\cA
-\cV$ denote the pushout of $\widehat\cA\times_{\cT_X}\cV$ by the
difference map $\Omega^1_X\times\Omega^1_X\xrightarrow{-}\Omega^1_X$.
Thus, a section of $\widehat\cA -\cV$ is represented by a pair
$(a,v)$ with $a\in\widehat\cA$, $v\in\cV$ satisfying
$\pi(a)=\pi(v)\in\cT_X$. Two pairs as above are equivalent if
their (componentwise) difference is of the form $(\alpha,\alpha)$
for some $\alpha\in\Omega^1_X$.

For $a\in\widehat\cA$, $v\in\cV$ with $\pi(a)=\pi(v)$, $f\in\cO_X$
let
\begin{equation}\label{formulas:vert-mult-diff}
f*(a,v) = (f*a, f*v),\ \ \ \partial(f) = \partial_{\widehat\cA}(f)
- \partial_{\cV}(f)\ .
\end{equation}

For $a_i\in\widehat\cA$, $v_i\in\cV$ with $\pi(a_i)=\pi(v_i)$ let
\begin{equation}\label{formulas:vert-bracket-diff}
[(a_1,v_1),(a_2,v_2)] = ([a_1,a_2],[v_1,v_2]),\ \ \
\langle(a_1,v_1),(a_2,v_2)\rangle = \langle a_1,a_2\rangle -
\langle v_1,v_2\rangle
\end{equation}
These operations are easily seen to descend to $\widehat\cA - \cV$.

The two maps $\Omega^1_X\to\widehat\cA - \cV$ given by
$\alpha\mapsto (i(\alpha),0)$ and $\alpha\mapsto (0,-i(\alpha))$
coincide; we denote their common value by
\begin{equation}\label{map:der-diff}
i : \Omega^1_X\to\widehat\cA - \cV \ .
\end{equation}

\begin{lemma}\label{lemma:comp-ECA-tors}
The formulas \eqref{formulas:vert-mult-diff}, \eqref{formulas:vert-bracket-diff}
together with \eqref{map:der-diff} determine a structure of a Courant
extension of $\cA$ on $\widehat\cA - \cV$. Moreover, the map
$\mathfrak{g}(\widehat\cA)\to\widehat\cA-\cV$ defined by $a\mapsto(a,0)$
induces an isomorphism $\mathfrak{g}(\widehat\cA-\cV)\cong\mathfrak{g}(\widehat\cA)$
of Courant extensions of ${\mathfrak{g}}(\cA)$ (by $\Omega^1_X$).
\end{lemma}
\begin{proof}
Left to the reader.
\end{proof}

\begin{proposition}\label{prop:CEXToppEVA}
Suppose that $\cA$ admits a vertex extension $\widehat\cA$; let $\ip$ denote the
induced invariant symmetric pairing on ${\mathfrak{g}}(\cA)$.
Then, the assignment $\cV\mapsto\widehat\cA - \cV$ extends to a functor
\begin{equation}\label{functor:EVAtoCEXT}
\widehat\cA-(\bullet) : \EVA_{\cO_X}\to\CEXT_{\cO_X}(\cA)_\ip
\end{equation}
which anti-commutes with the respective actions of $\ECA_{\cO_X}$ on
$\EVA_{\cO_X}$ and $\CEXT_{\cO_X}(\cA)_\ip$. In particular, the functor
\eqref{functor:EVAtoCEXT} is an equivalence of stacks in groupoids.
The isomorphism classes of the $\ECA_{\cO_X}$-torsors  $\EVA_{\cO_X}$ and
$\CEXT_{\cO_X}(\cA)_\ip$ are opposite as elements of
$H^2(X;\Omega^2_X\to\Omega^{3,cl}_X)$.
\end{proposition}
\begin{proof}
This is clear from the construction of \ref{ss:comp-ECA-tors} and
Lemma \ref{lemma:comp-ECA-tors}.
\end{proof}


\begin{theorem}\label{main-thm}
The class of $\EVA_{\cO_X}$ in $H^2(X;\Omega^2_X\to\Omega^{3,cl}_X)$
is equal to $\ch_2(\Omega^1_X)$.
\end{theorem}
\begin{proof}
According to Proposition \ref{prop:the-algd}, $\cA_{\Omega^1_X}$ (the Atiyah algebra of $\Omega^1_X$)
admits a canonical vertex extension with the induced symmetric pairing on the Lie algebra
$\shEnd_{\cO_X}(\Omega^1_X) =
{\mathfrak{g}}(\cA_{\Omega^1_X})$ canonical invariant symmetric pairing
\footnote{see \ref{ssection:vb} for terminology and notations}
given by the trace of the product of endomorphisms. By Corollary \ref{cor:main-for-vb},
the class of $\CEXT_{\cO_X}(\cA_{\Omega^1_X})_{\Tr}$ is equal
to $-\ch_2(\Omega^1_X)$. The claim follows from Proposition \ref{prop:CEXToppEVA}.
\end{proof}

\section{Algebroids over the de Rham complex}
All of the notions of the preceding sections generalize in an
obvious way to differential graded manifolds (i.e. manifolds
whose structure sheaves are a sheaves of commutative differential graded
algebras).

For a manifold $X$ let $X^\sharp$ denote the differential graded
manifold with the underlying space $X$ and the structure sheaf $\cO_{X^\sharp}$
the de Rham complex $\DR$. In other words, $\cO_{X^\sharp} = \bigoplus_i\Omega^i_X[-i]$
(as a sheaf of graded algebras). We will denote by $\partial_{\cO_{X^\sharp}}$
the derivation given by the de Rham differential.

\subsection{The structure of $\cT_{X^\sharp}$}
The tangent sheaf of $X^\sharp$ (of derivations of $\cO_{X^\sharp}$),
$\cT_{X^\sharp}$ is a sheaf of differential graded Lie algebras with the
differential $\partial_{\cT_{X^\sharp}} = [\partial_{\cO_{X^\sharp}},\ \ ]$
(note that $\partial_{\cO_{X^\sharp}}\in\cT_{X^\sharp}^1$).

\subsubsection{}
Let $\widetilde\cT_X$ denote the cone of the identity endomorphism
of $\cT_X$. That is, $\widetilde{\cT_X}^i=\cT_X$ for $i=-1,0$ and
zero otherwise. The only nontrivial differential is the identity map.
The complex $\widetilde\cT_X$ has the canonical structure of a
sheaf of differential graded Lie algebras.

The natural action of $\cT_X$ (respectively $\cT_X[1]$) on $\cO_{X^\sharp}$
by the Lie derivative (respectively by interior product) gives rise
to the injective map of
DGLA
\begin{equation}\label{action-tau}
\tau : \widetilde{\cT_X}\to\cT_{X^\sharp} \ .
\end{equation}

The action $\tau$ extends in the canonical way to a structure
of a Lie $\cO_{X^\sharp}$-algebroid on $\cO_{X^\sharp}\otimes_\mathbb{C}\widetilde\cT_X$
with the anchor map
\begin{equation}\label{tau-DR}
\tau_{\cO_{X^\sharp}} : \cO_{X^\sharp}\otimes_\mathbb{C}\widetilde\cT_X
\to\cT_{X^\sharp}
\end{equation}
the canonical extension $\tau$ to a $\cO_{X^\sharp}$-linear map. Note that
$\tau_{\cO_{X^\sharp}}$ is surjective, i.e. the Lie $\cO_{X^\sharp}$-algebroid
$\cO_{X^\sharp}\otimes_\mathbb{C}\widetilde\cT_X$ is transitive. We denote this
algebroid by $\widetilde\cT_{X^\sharp}$.

Let $\cT_{X^\sharp/X}\subset\cT_{X^\sharp}$ denote the normalizer of
$\cO_X\subset\cO_{X^\sharp}$. Since the action of $\cT_X[1]$ is
$\cO_X$-linear, the map $\tau$ restricts to
\[
\tau : \widetilde{\cT_X}^{-1} = \cT_X[1]\to\cT_{X^\sharp/X}
\]
and (the restriction of) $\tau_{\cO_{X^\sharp}}$ factors through the map
\begin{equation}\label{tau-sub}
\cO_{X^\sharp}\otimes_{\cO_X}\cT_X[1]\to\cT_{X^\sharp/X}
\end{equation}
which is easily seen to be an isomorphism.

Since the action $\tau$ is $\cO_X$-linear modulo
$\cT_{X^\sharp/X}$, $\tau_{\cO_{X^\sharp}}$ induces the map
\begin{equation}\label{tau-fac}
\cO_{X^\sharp}\otimes_{\cO_X}\cT_X\to\cT_{X^\sharp}/\cT_{X^\sharp/X}
\end{equation}
which is easily seen to be an isomorphism.

Therefore, there is an exact sequence of graded
$\cO_{X^\sharp}$-modules
\begin{equation}\label{ses-der}
0\to\cO_{X^\sharp}\otimes_{\cO_X}\cT_X[1]\to\cT_{X^\sharp}\to
\cO_{X^\sharp}\otimes_{\cO_X}\cT_X\to 0\ .
\end{equation}
The composition
\[
\cO_{X^\sharp}\otimes_{\cO_X}\cT_X[1]\to\cT_{X^\sharp}
\xrightarrow{\partial_{\cT_{X^\sharp}}}
\cT_{X^\sharp}[1]\to\cO_{X^\sharp}\otimes_{\cO_X}\cT_X[1]
\]
is the identity map.

The natural action of $\cT_{X^\sharp}$ on $\cO_{X^\sharp} = \Omega^\bullet_X$ restricts
to the action of $\cT_{X^\sharp}{}^0$ on $\cO_X$ and $\Omega^1_X$. The action of
$\cT_{X^\sharp}{}^0$ on $\cO_X$ gives rise to the map
$\cT_{X^\sharp}{}^0\to\cT_X$ which, together with the natural Lie bracket on
$\cT_{X^\sharp}{}^0$, endows the latter with a structure of a Lie $\cO_X$-algebroid.

The action of $\cT_{X^\sharp}{}^0$ on $\Omega^1_X$ gives rise to the map
\begin{equation}\label{map:TtoA}
\cT_{X^\sharp}{}^0\to\cA_{\Omega^1_X} \ ,
\end{equation}
where $\cA_{\Omega^1_X}$ denotes the Atiyah algebra of $\Omega^1_X$.

\begin{lemma}
The map \eqref{map:TtoA} is an isomorphism of Lie $\cO_X$-algebroids.
\end{lemma}

\subsection{Exact Courant $\cO_{X^\sharp}$-algebroids}

\begin{proposition}\label{prop:ECA-DR-triv}
Every exact Courant $\cO_{X^\sharp}$-algebroid admits a unique flat
connection.
\end{proposition}
\begin{proof}
Consider an exact Courant $\cO_{X^\sharp}$-algebroid
\[
0\to\Omega^1_{X^\sharp}\xrightarrow{i}\cQ\xrightarrow{\pi}\cT_{X^\sharp}\to 0
\]
Note that, since $\Omega^1_{X^\sharp}$ is concentrated in non-negative
degrees, the map $\pi : \cQ^{-1}\to\cT_{X^\sharp}{}^{-1}$ is
an isomorphism. Since $\cT_{X^\sharp}{}^{-1}$ generates $\cT_{X^\sharp}$ as a
DG-module over $\cO_{X^\sharp}$, the splitting is unique if it exists.

To establish the existence it is necessary and sufficient to show
that the restriction of the anchor map to the DG-submodule of
$\cQ$ generated by $\cQ^{-1}$ is an isomorphism.

Note that the map $\tau : \widetilde\cT_X\to\cT_{X^\sharp}$ lifts in a
unique way to a morphism of complexes $\widetilde\tau :
\widetilde\cT_X\to\cQ$. The map $\widetilde\tau$ is easily seen to be a morphism of
DGLA. Let $\cQ'$ denote the $\cO_{X^\sharp}$-submodule of $\cQ$ generated by
the image of $\widetilde\tau$ (i.e. the DG $\cO_{X^\sharp}$-submodule
generated by $\cQ^{-1}$).

Since
\[
\widetilde\tau :\cT_X[1] = \widetilde\cT_X^{-1}\to\cQ'
\]
is $\cO_X$-linear it extends to the map
\begin{equation}\label{Q-sub}
\cO_{X^\sharp}\otimes_{\cO_X}\cT_X[1]\to\cQ'
\end{equation}
such that the composition
\[
\cO_{X^\sharp}\otimes_{\cO_X}\cT_X[1]\to\cQ'\xrightarrow{\pi}\cT_{X^\sharp}
\]
coincides with the composition of the isomorphism
\eqref{tau-sub}with the inclusion into $\cT_{X^\sharp}$. Therefore,
\eqref{Q-sub} is a monomorphism whose image will be denoted
$\cQ''$, and $\pi$ restricts to an isomorphism of $\cQ''$ onto
$\cT_{X^\sharp/X}$.

Since
\[
\widetilde\tau : \cT_X = \widetilde\cT_X^0\to\cQ'\to\cQ'/\cQ''
\]
is $\cO_X$-linear it extends to the map
\begin{equation}\label{Q-fac}
\cO_{X^\sharp}\otimes_{\cO_X}\cT_X\to\cQ'/\cQ''
\end{equation}
which is surjective (since $\cQ'/\cQ''$ is generated as a
$\cO_{X^\sharp}$-module by the image of $\widetilde\cT_X^0$ under $\widetilde\tau$),
and such that the composition
\[
\cO_{X^\sharp}\otimes_{\cO_X}\cT_X\to\cQ'/\cQ''\xrightarrow{\pi}
\cT_{X^\sharp}/\cT_{X^\sharp/X}
\]
coincides with the isomorphism \eqref{tau-fac}. Therefore,
\eqref{Q-fac} is an isomorphism.

Now, the exact sequence \eqref{ses-der} implies that $\pi$
restricts to an isomorphism $\cQ'\cong\cT_{X^\sharp}$. The desired
splitting is the inverse isomorphism. It is obviously compatible
with brackets, hence, is a flat connection.
\end{proof}

\begin{corollary}\label{ECA-DR-final}
$\ECA_{\cO_{X^\sharp}}$ is equivalent to the (final) stack
$X\supseteq U\mapsto [0]$, where $[0]$ is the category with one
object and one morphism.
\end{corollary}

\begin{corollary}\label{corollary:existence-uniqueness}
An exact vertex $\cO_{X^\sharp}$-algebroid exists and is unique up to
canonical isomorphism.
\end{corollary}
\begin{proof}
Since $\EVA_{\cO_{X^\sharp}}$ is an affine space under $\ECA_{\cO_{X^\sharp}}$
the uniqueness (local and global) follows from Corollary
\ref{ECA-DR-final}. Local existence and uniqueness implies
global existence.
\end{proof}

\subsection{The canonical vertex $\cO_X$-algebroid}
Let $\cV^\bullet$ denote the unique exact vertex $\cO_{X^\sharp}$-algebroid.
The degree zero component of the exact sequence
\[
0\to\Omega^1_{X^\sharp}\to\cV^\bullet\to\cT_{X^\sharp}\to0
\]
is canonically isomorphic to
\[
0\to\Omega^1_X\to\cV^0\to\cA_{\Omega^1_X}\to 0
\]
using the canonical isomorphisms $\Omega^1_{X^\sharp}{}^0\isomo\Omega^1_X$
and \eqref{map:TtoA}.

\begin{proposition}\label{prop:the-algd}
The vertex $\cO_{X^\sharp}$-algebroid structure on $\cV$ restricts to a
structure of a vertex extension of $\cA_{\Omega^1_X}$ on $\cV^0$.
The induced symmetric pairing on $\shEnd_{\cO_X}(\Omega^1_X)$ is given by
the trace of the product of endomorphisms.
\end{proposition}
\begin{proof}
The first statement is left to the reader.

The degree $-1$ component of the anchor map $\cV^\bullet\to\cT_{X^\sharp}$
is an isomorphism whose inverse gives the map $\cT_X[1]\to\cV^\bullet$.
Combined with the $*$-multiplication by $\cO_{X^\sharp}^1 = \Omega^1_X$ it
gives the map $\Omega^1_X\otimes_\mathbb{C}\cT_X\to\cV^0$ and the commutative
diagram
\[
\begin{CD}
\Omega^1_X\otimes_\mathbb{C}\cT_X @>{*}>> \cV^0 \\
@VVV @VVV \\
\shEnd_{\cO_X}(\Omega^1_X) @>>> \cA_{\Omega^1_X}
\end{CD}
\]
where the left vertical map is the canonical one.

For $\phi,\psi\in\shEnd_{\cO_X}(\Omega^1_X)$ represented, respectively,
by $\alpha\otimes\xi$ and $\beta\otimes\eta$ with $\alpha,\beta\in\Omega^1_X[-1]$
and $\xi,\eta\in\cT_X[1]$, the pairing $\langle\phi,\psi\rangle$ is calculated
using \eqref{pairing}:
\begin{multline*}
\langle\phi,\psi\rangle = \langle\alpha *\xi,\beta *\eta\rangle_\cV =
\alpha\langle\xi,\beta *\eta\rangle +\pi(\xi)(\pi(\beta *\eta)(\alpha)) = \\
\pi(\xi)(\beta\pi(\eta)(\alpha)) =
\pi(\xi)(\beta)\pi(\eta)(\alpha) -
\beta\pi(\xi)(\pi(\eta)(\alpha)) = \iota_\xi\beta\cdot\iota_\eta\alpha =
\Tr(\phi\psi)
\end{multline*}
(with $\langle\xi,\beta *\eta\rangle = 0$ and $\pi(\xi)(\pi(\eta)(\alpha)) =0$
since both have negative degrees).
\end{proof}

\appendix
\section{Characteristic classes of Lie algebroids}\label{section:char-classes}
Below, we recall the basic definitions and facts regarding Lie algebroids as well
as the construction of the higher Chern-Simons forms (\cite{CS}). As an example we calculate explicitly the \v{C}ech-de Rham representative of the first Pontryagin class.

\subsection{Lie algebroids}\label{ssection:lie-algebroids}
We refer the reader to \cite{M} for further details on Lie algebroids in the differential-geometric context.

\subsubsection{Definitions}
Suppose that $X$ is a manifold. A \textit{Lie algebroid} on $X$ is a sheaf $\cA$ of $\cO_X$ modules equipped with
an $\cO_X$-linear map $\pi: \cA\to\cT_X$ called the \textit{anchor map}, and a $\mathbb{C}$-linear pairing
$[\ ,\ ]:\cA\otimes_{\mathbb{C}}\cA\to\cA$ such that
\begin{enumerate}
\item the pairing $[\ ,\ ]$ is a Lie bracket (i.e. it is skew-symmetric and satisfies
the Jacobi identity),
\item the map $\pi$ is a morphism of Lie algebras (i.e. it commutes with the respective brackets),
\item the Leibniz rule holds, i.e. for $a_1,a_2\in\cA$, $f\in\cO_X$ the identity
$[a_1,fa_2] = \pi(a_1)(f)a_2 + f[a_1,a_2]$ is satisfied.
\end{enumerate}
We denote by $\mathfrak{g}(\cA)$ the kernel of the anchor map $\pi$. This is an $\cO_X$-Lie algebra.

A morphism of Lie algebroids is an $\cO_X$-linear map of Lie algebras which commutes
with the respective anchor maps.

A Lie algebroid is called \textit{transitive} if the anchor map is surjective.

\subsubsection{Connections and curvature}\label{connections for Lie}
A \textit{connection} on a (transitive) Lie algebroid
$\cA$ is an $\cO_X$-linear splitting $\nabla : \cT_X\to\cA$ of the anchor map (i.e.
$\pi\circ\nabla=\id$). A  connection is called \textit{flat} if it commutes with
the respective brackets (i.e. if it is a morphism of Lie algebras).

For a connection $\nabla$ the formula $c(\nabla)(\xi_1,\xi_2) := [\nabla(\xi_1),\nabla(\xi_2)]-\nabla([\xi_1,\xi_2])$ defines
the $\cO_X$-linear map $c(\nabla) : \bigwedge^2\cT\to\mathfrak{g}(\cA)$ called
the \textit{curvature} of the connection $\nabla$. Clearly, $\nabla$ is flat if and
only if $c(\nabla)=0$.

\subsubsection{Examples}
The tangent sheaf $\cT_X$ (the anchor map being the identity) is the final object in
in the category of Lie algebroids on $X$.

An algebroid with the trivial anchor map is the same thing as a sheaf of $\cO_X$-Lie
algebras.

An action of a Lie algebra
$\mathfrak{g}$ on $X$, the action given by the morphism of Lie algebras $\alpha:\mathfrak{g}\to\Gamma(X;\cT_X)$ gives rise to
a structure of a Lie algebroid on $\cO_X\otimes\mathfrak{g}$ in a natural way: the anchor map is given by
$f\otimes a\mapsto f\alpha(a)$ and the bracket is defined by $[f_1\otimes a_1,f_2\otimes a_2] =
f_1 f_2\otimes [a_1,a_2] + f_1\alpha(a_1)(f_2)\otimes a_2 - f_2\alpha(a_2)(f_1)\otimes a_1$.

\subsubsection{Atiyah algebras}
An important class of examples of transitive Lie algebroids is comprised of Atiyah algebras. Suppose that $G$ is a Lie group
with Lie algebra $\mathfrak{g}$ and $p: P\to X$ is a principal $G$ bundle. The \textit{Atiyah algebra} of $P$, denoted
$\cA_P$, is the sheaf whose local sections are pairs $(\tilde\xi,\xi)$, where $\xi$ is a (locally defined) vector field on
$X$ and $\tilde\xi$ is a $G$-invariant vector field on $P$ which lifts $\xi$. Thus, $\cA_P=\left(p_*\cT_P\right)^G$
with the induced bracket and the anchor map is given by
\[
\cA_{\cE}=\left(p_*\cT_{\cE}\right)^G\stackrel{dp}
{\longrightarrow}\left(p_*p^*\cT_X\right)^G=\cT_X \ .
\]
In terms of the preceding description the anchor is give by the projection $(\tilde\xi,\xi)\mapsto\xi$.
In this case $\mathfrak{g}(\cA)$ is the $\cE$ twist of $\cO_X\otimes\mathfrak{g}$.

The Atiyah algebra of a vector bundle, which is just the Atiyah algebra of the corresponding $GL$-bundle, admits
another description. Namely, for a vector bundle $\cF$ let $\Diff^{\leq n}(\cF,\cF)$ denote the sheaf of differential
operators of order $n$ acting on $\cF$. Then $\cA_{\cF}$ is determined by the following pull-back diagram
\[
\begin{CD}
0 @>>> \shEnd_{\cO_X}(\cF) @>>> \cA_{\cF} @>\pi>> \cT_X @>>> 0 \\
& & @VVV @VVV @VV{\id\otimes\vac}V \\
0 @>>>\Diff^{\leq 0}(\cF,\cF) @>>> \Diff^{\leq 1}(\cF,\cF) @>>> \cT_X\otimes\shEnd_{\cO_X}(\cF)
@>>> 0
\end{CD}
\]
as the sheaf of differential operators of order one with scalar principal symbol.

\subsubsection{Pull-back of algebroids}
Recall that, for a Lie algebroid $\cA$ on $X$ and a map $\phi : Y\to X$, the pull-back
$\phi^+\cA$ is defined by the Cartesian diagram
\[
\begin{CD}
\phi^+\cA @>>> \phi^*\cA \\
@VVV @VV{\id\otimes\pi}V \\
\cT_Y @>{d\phi}>> \phi^*\cT_X
\end{CD}
\]
Thus, a section of $\phi^+\cA$ is a pair $(f\otimes a,\xi)$, where
$f\in\cO_Y$, $a\in\cA$, $\xi\in\cT_Y$, so that $f\otimes a\in\phi^*\cA =
\cO_Y\otimes_{\phi^{-1}\cO_X}\cA$, and $f\otimes\pi(a) = d\phi(\xi)$.
The anchor map is given by the projection $(f\otimes a,\xi)\mapsto\xi$.

The bracket on $\phi^+\cA$ is the unique one which restricts to the bracket
on $\phi^{-1}\cA$ (induced by that on $\cA$) and obeys the Leibniz rule.
Note that $\mathfrak{g}(\phi^+\cA)\isomo \phi^*\mathfrak{g}(\cA)$ as sheaves of
$\cO_Y$-Lie algebras.

\subsection{Higher Chern-Simons forms}\label{ssec:CS}
Suppose that $\cU$ is a cover of $X$ by open subsets. Let
$X_0 = \coprod_{U\in\cU} U$, $\epsilon : X_0\to X$ the map induced by the
inclusions $i_U : U\hookrightarrow X$. Let $X_i = X_0\times_X\dots\times_X X_0$
denote the $(i+1)$-fold product which will be indexed by $\{0, \ldots, i\}$.
For $j = 0,\ldots,i$ let $\pr_j : X_i\to X_0$ (respectively, $s_j : X_i\to X_{i-1}$)
denote the projection onto (respectively, along) the $j^{\text th}$ factor.

The collection of all $X_i$ together with
the projections $s_j$ and the diagonal maps is a simplicial manifold denoted
$X_\bullet$, and $\epsilon$ extends to the map $\epsilon : X_\bullet\to X$,
where the latter is regarded as a constant simplicial object.

A sheaf $F$ on $X$ gives rise to a simplicial sheaf $\epsilon^*F$
on $X_\bullet$ so that the complex (associated to the simplicial abelian group)
$\Gamma(X_\bullet;F)$ is the complex $\check C(\cU;F)$ of $F$-valued cochains
on the cover $\cU$. We will denote by $F_i$ the restriction of $\epsilon^*F$
to $X_i$. Note that there is a canonical isomorphism $F_i\isomo\pr_j^*F_0$.

Let $\Delta^i$ denote the standard $i$-dimensional simplex:
\[
\Delta^i = \lbrace \vec{t}=(t_0,\ldots,t_i)\vert\ t_j\in\mathbb{R},
\ \sum_{j=0}^i t_j\leq 1\rbrace \ .
\]
Integration over $\Delta^i$ gives rise to the map
\[
\int_{\Delta_i} : \Gamma(X_i\times\Delta^i;\Omega^a_{X_i\times\Delta^i})
\to\Gamma(X_i; \Omega^{a-i}_{X_i})
\]
which satisfies the Stokes formula
\[
\int_{\Delta_i}d\alpha = d\int_{\Delta_i}\alpha +
\sum_{j=0}^i (-1)^j\int_{\Delta^{i-1}}(\id\times\partial_j)^*\alpha
\]
where $\partial_j : \Delta^{i-1}\to\Delta_i$ is the map
$(t_0,\ldots,t_{i-1})\mapsto(t_0,\ldots,t_{j-1},0,t_j,\ldots,t_{i-1})$.

The forms on the product $X_i\times\Delta^i$ are naturally bi-graded:
$\Omega^a_{X_i\times\Delta^i}=\Omega^{a-i,i}_{X_i\times\Delta^i}+
\Omega^{a-i+1,i-1}_{X_i\times\Delta^i}+\dots+\Omega^{a,0}_{X_i\times\Delta^i}$.
The ``integration over $\Delta^i$'' map is non-trivial only on the first summand
(of maximal $\Delta^i$-degree) which, on the other hand lies in the kernel of
the restriction maps $\partial_j^*$. Hence, for
$\alpha = \sum_{k=0}^i\alpha^{a-k,k}$ with
$\alpha^{a-k,k}\in\Gamma(X_i\times\Delta^i;\Omega^{a-k,k}_{X_i\times\Delta^i})$
the Stokes formula takes the form
\begin{equation}\label{formula:Stokes}
\int_{\Delta_i}d\alpha = d\int_{\Delta_i}\alpha^{a-i,i} +
\sum_{j=0}^i (-1)^j\int_{\Delta^{i-1}}(\id\times\partial_j)^*\alpha^{a-i+1,i-1}
\end{equation}

Suppose that $\cA$ is a Lie algebroid on $X$, locally free of finite rank over
$\cO_X$, and let $\mathfrak{g} := \mathfrak{g}(\cA)$. Let
$\widetilde\cA_i = \pr_{X_i}^+\cA_i$, where $\pr_{X_i}:X_i\times\Delta^i\to X_i$
denotes the projection, so that
\[
\widetilde\cA_i =
\left(\cA_i\boxtimes\cO_{\Delta^i}\oplus\cO_{X_i}\boxtimes\cT_{\Delta^i}\right)
\otimes_{\cO_{X_i}\boxtimes\cO_{\Delta^i}}\cO_{X_i\times\Delta^i}
\]
as a $\cO_{X_i\times\Delta^i}$-module with the obvious Lie algebroid structure.

For a connection $\nabla_0$ on $\cA_0$ let $\widetilde\nabla_i$ denote
the connection on $\widetilde\cA_i$ given by the formula
\[
\widetilde\nabla_i = \sum_{j=0}^i \pr_j^*\nabla_0\boxtimes t_j +
1\boxtimes\id \ .
\]

The curvature of $\widetilde\nabla_i$, $c(\widetilde\nabla_i)\in
\Omega^2_{X_i\times\Delta^i}\otimes\widetilde{\mathfrak{g}}_i$ decomposes
as $c(\widetilde\nabla_i) = c(\widetilde\nabla_i)^{2,0} +
c(\widetilde\nabla_i)^{1,1} + c(\widetilde\nabla_i)^{0,2}$ with
$c(\widetilde\nabla_i)^{2-k,k}\in
\Omega^{2-k,k}_{X_i\times\Delta^i}\otimes\widetilde{\mathfrak{g}}_i$.
Since $\widetilde\nabla_i$ is flat along $\Delta^i$, the $(0,2)$-component
of the curvature form vanishes.

Suppose that $P : \mathfrak{g}^{\otimes_{\cO_X} d}\to\cO_X$
is an $\cA$-invariant map. Then, the map
$\widetilde{P}_i:\widetilde{\mathfrak{g}}_i^{\otimes_{\cO_{X_i\times\Delta^i}} d}\to
\cO_{X_i\times\Delta^i}$ is $\widetilde\cA_i$-invariant. The invariance of
$\widetilde{P}_i$ implies that the differential form
$\widetilde{P}_i(c(\widetilde\nabla_i)^{\wedge d})
\in\Gamma(X_i\times\Delta^i;\Omega^{2d}_{X_i\times\Delta^i})$ is closed.
In view of the remarks above,
\begin{multline*}
\widetilde{P}_i(c(\widetilde\nabla_i)^{\wedge d}) =
\widetilde{P}_i((c(\widetilde\nabla_i)^{2,0} +
c(\widetilde\nabla_i)^{1,1})^{\wedge d}) = \\
\sum_{k=0}^d\binom{d}{k}\widetilde{P}_i((c(\widetilde\nabla_i)^{2,0})^{\wedge(d-k)}
\wedge(c(\widetilde\nabla_i)^{1,1})^{\wedge k})
\end{multline*}
with the $k^{\text th}$ summand, denoted $\widetilde{P}_i(c(\widetilde\nabla_i)^{\wedge d})^{2d-k,k}$ below, homogeneous of bi-degree $(2d-k,k)$.

Recall that $s_j : X_i\to X_{i-1}$ denotes the projection along the $j^{\text{th}}$ factor.
The Lie algebroids $(\id\times\partial_j)^+\widetilde\cA_i$ and
$(s_j\times\id)^+\widetilde\cA_{i-1}$ are canonically isomorphic in a way
compatible with the connections (induced by)
$\widetilde\nabla_i$ on the former and $\widetilde\nabla_{i-1}$ on the latter.
It follows that
\[
(\id\times\partial_j)^*\widetilde{P}_i(c(\widetilde\nabla_i)^{\wedge d})^{2d-k,k}
=(s_j\times\id)^*\widetilde{P}_{i-1}(c(\widetilde\nabla_{i-1})^{\wedge d})^{2d-k,k} \ .
\]

Using
\begin{enumerate}
\item
\begin{multline*}
\int_{\Delta^{i-1}}(\id\times\partial_j)^*
\widetilde{P}_i(c(\widetilde\nabla_i)^{\wedge d})^{2d-i+1,i-1} = \\
\int_{\Delta^{i-1}}(s_j\times\id)^*
\widetilde{P}_{i-1}(c(\widetilde\nabla_{i-1})^{\wedge d})^{2d-i+1,i-1} = \\
s_j^*\int_{\Delta^{i-1}}
\widetilde{P}_{i-1}(c(\widetilde\nabla_{i-1})^{\wedge d})^{2d-i+1,i-1} \ ,
\end{multline*}

\item the fact that the differential in the \v{C}ech complex is given by
$\check\partial = \sum_j (-1)^j s_j^*$ \ ,

\item the Stokes formula \eqref{formula:Stokes} and

\item noticing that the left hand side of the latter vanishes because
$\widetilde{P}_i(c(\widetilde\nabla_i)^{\wedge d})$ is closed
\end{enumerate}
one sees that
\[
\sum_i\int_{\Delta^i}\widetilde{P}_i(c(\widetilde\nabla_i)^{\wedge d})^{d-i,i}
\in\check{C}^\bullet(\cU;\tau_{\leq 2d}\Omega^{\geq d}_X)
\]
is a cocycle of total degree $2d$ in the \v{C}ech--de Rham complex which extends $P(c(\nabla_0)^{\wedge d})\in\check{C}^0(\cU;\Omega^{2d,cl}_X)$, and whose class in
$H^{2d}(X;\tau_{\leq 2d}\Omega^{\geq d}_X)$ is independent of the choices
(of the covering $\cU$ and the connection $\nabla_0$) made. (In fact, independence
of the choice of the connection follows from the above construction applied to
the simplicial manifold $X_\bullet$ in place of $X$, the algebroid
$\cA_\bullet = \epsilon^*\cA$ and the connection $\pr_\bullet^*\nabla$ induced
by a choice of a connection $\nabla$ on $\cA_0$.)

\subsection{The first Pontryagin class}\label{ssection:pervi-pont}
Below we will carry out the above calculation in the situation when
$P = \ip$ is a $\cA$-invariant symmetric pairing on $\mathfrak g$. The
corresponding class in $H^4(X;\tau_{\leq 4}\Omega^{\geq 2}_X)$ is called
the (first) Pontryagin class of the pair $(\cA,\ip)$ and will be denoted
$\pont(\cA,\ip)$. It follows from the discussion in \ref{ssec:CS} that
$\pont(\cA,\ip)$ is represented, in terms of a covering $\cU$ and a connection
$\nabla_0$ by the cocycle $\pont^{4,0}+\pont^{3,1}+\pont^{2,2}$ with
$\pont^{i,4-i}\in\check{C}^{4-i}(\cU;\Omega^i_X)$ (where we have suppressed
the dependence on the data).

Clearly, $\pont^{4,0}=\langle c(\nabla_0)\wedge c(\nabla_0)\rangle$.

\subsubsection{Calculation of $\pont^{3,1}$}
Let $A = \pr_1^*\nabla_0 - \pr_0^*\nabla_0$, $\nabla := \pr_0^*\nabla_0$,
$t := t_1$. With these notations $\widetilde\nabla_1 =
\nabla + At + d_t$,
\[
c(\widetilde\nabla_1) = \frac12[\widetilde\nabla_1,\widetilde\nabla_1]
=\frac12[\nabla,\nabla] + [\nabla,A]t + \frac12[A,A]t^2 - Adt \ ,
\]
so that
\[
(c(\widetilde\nabla_1)\wedge c(\widetilde\nabla_1))^{3,1} =
-([\nabla,\nabla] + 2[\nabla,A]t + [A,A]t^2)\wedge Adt
\]
which, upon integration, gives
\begin{multline*}
\pont^{3,1} = -(\langle[\nabla,\nabla]\wedge A\rangle +
\langle[\nabla,A]\wedge A\rangle + \frac13\langle[A,A]\wedge A\rangle) \\
=
-2(\langle c(\nabla)\wedge A\rangle + \frac12\langle[\nabla,A]\wedge A\rangle
+ \frac16\langle[A,A]\wedge A\rangle)\ .
\end{multline*}

\subsubsection{Calculation of $\pont^{2,2}$}
Let $A_{ij} = \pr_j^*\nabla_0 - \pr_i^*\nabla_0$, $\nabla := \pr^*_1\nabla_0$.
With these notations
\begin{multline*}
\widetilde\nabla_2 =
\sum_{j=0}^2 \pr_j^*\nabla_0\cdot t_j + d_t =
\sum_{j=0}^2 (\pr_1^*\nabla_0 + A_{j1})\cdot t_j + d_t = \\
\nabla + A_{10}t_0 + A_{12}t_2 + d_t =
\nabla - A_{01}t_0 + A_{12}t_2 + d_t \ ,
\end{multline*}
\[
c(\widetilde\nabla_2) = \frac12[\widetilde\nabla_2,\widetilde\nabla_2] =
\dots + A_{01}dt_0 - A_{12}dt_2)
\]
(ignoring the terms which do not contain $dt_1$ and $dt_2$), so that
\[
(c(\widetilde\nabla_2)\wedge c(\widetilde\nabla_2))^{2,2} =
2A_{01}\wedge A_{12}\wedge dt_0\wedge dt_2
\]
which, upon integration, gives
\[
\pont^{2,2} = -\langle A_{01}\wedge A_{12}\rangle \ .
\]

\subsection{Pontryagin class for vector bundles}\label{ssection:vb}
Suppose that $\cE$ is a $GL_n$-torsor on $X$. Let $\cA_{\cE}$ denote the Atiyah algebra of
$\cE$. Thus, $\cA_\cE$ is a transitive Lie algebroid with $\mathfrak{g}(\cA_\cE) =
{\mathfrak{gl}}_n^\cE$. The Lie algebra $\mathfrak{gl}_n$ carried the canonical invariant symmetric
pairing given by the trace of the product of endomorphisms which induces the canonical
$\cA_\cE$-invariant symmetric pairing on $\mathfrak{gl}_n^\cE$. We will denote this pairing by $\ip_{can}$ and write $\pont(\cE)$ for $\pont(\cA_{\cE},\ip_{can})$.

Recall that the Chern character form of a connection $\nabla$ on $\cE$,  denoted
$\ch(\cE,\nabla) = \sum_n\ch_n(\cE,\nabla)$, where $\ch_n(\cE,\nabla)$ a form of degree $2n$, is defined by
\[
\ch(\cE,\nabla) = \Tr(\exp(c(\nabla)) = \sum\frac1{n!}\Tr(c(\nabla)^{\wedge n})
\]
In particular, the characteristic class $\ch_2$ is  related to the Pontryagin class by
\begin{equation}\label{formula:pont-ch}
\pont(\cE) = 2\ch_2(\cE) \ .
\end{equation}


\begin{thebibliography}{ABCDE}
\bibitem[Bl81]{B} S.~Bloch, The dilogarithm and extensions of Lie algebras, in \textit{Algebraic
$K$-theory, Evanston 1980 (Proc. Conf., Northwestern Univ., Evanston, Ill. ,1980)},
Lecture Notes in Math., 854, Springer, 1981, pp. 1--23.
\bibitem[Br94]{Br} L.~Breen, On the classification of 2-gerbes and 2-stacks, Asterisque \textbf{225}, (1994).
\bibitem[BD04]{BD} A.~Beilinson, V.~Drinfeld, \textit{Chiral algebras}, A.M.S. Colloquium Publications \textbf{51}, 2004.
\bibitem[C90]{C} T.J.~Courant, Dirac manifolds, {\em Trans. A.M.S.} \textbf{319} (1990), 631--661.
\bibitem[CS74]{CS} S.--S.~Chern, J.~Simons, Characteristic forms and geometric
invariants, {\em Annals Math.} \textbf{99} (1974), 48--69.
\bibitem[GMS04]{GMS} V.~Gorbunov, F.~Malikov, V.~Schechtman, Gerbes of chiral
differential operators II, Vertex algebroids, \textit{Invent. Math.} \textbf{155} (2004), 605--680.
\bibitem[LWX97]{LWX} Zhang-Ju Liu, A.~Weinstein, Ping Xu, Manin triples for
Lie bialgebroids, {\em J. Differential Geom.} \textbf{45} (1997), 547--574.
\bibitem[Mac87]{M} K.~Mackenzie, \textit{Lie groupoids and Lie algebroids in differential geometry},
London Mathematical Society Lecture Notes Series \textbf{124}, Cambridge University Press, 1987.
\bibitem[MSV99]{MSV} F.~Malikov, V.~Schechman, A.~Vaintrob, Chiral de Rham complex,
\textit{Comm. Math. Phys.} \textbf{204} (1999), 429--473.
\bibitem[\v Sev]{S} P.~\v Severa, letters to A.~Weinstein.
\end{thebibliography}
\end{document}